\documentclass[11pt
]{article}

\usepackage{amssymb}
\usepackage{amsfonts,amstext,amsmath,amsthm,
latexsym,mathrsfs,amsbsy,mparhack}
\usepackage{amssymb}
\usepackage{accents}

\usepackage{graphicx}
\usepackage[dvipsnames]{xcolor}

\usepackage[bbgreekl]{mathbbol}
\usepackage{bold-extra}

\usepackage{marginnote}
\usepackage{multicol}

\usepackage{makeidx}

\usepackage[
english,greek]{babel}

\usepackage[b]{esvect}

\usepackage{tocloft,relsize}

\DeclareSymbolFontAlphabet{\mathbbl}{bbold}

\usepackage{doi}

\usepackage{hyperref}

\usepackage{tikz}

\usepackage{flowchart}
\usetikzlibrary{arrows}

\usepackage{bm}

\usepackage{array,longtable}
\newcolumntype{C}{>{$}c<{$}}  
\newcolumntype{L}{>{$}l<{$}}  
\usepackage
{lineno}

\input{80f.sty}
\input{80fgeN.sty}

\hypersetup
{%
hypertex=true,
hypertexnames=true,
hidelinks,
pagebackref=true,
citecolor = Magenta,
colorlinks=true, 
}

\usepackage{cleveref}

\begin{document}

\selectlanguage{english}

\setcounter{tocdepth}{2}


\title{Notes on the equiconsistency of ZFC
without the Power Set axiom and  
second order PA%
%
\thanks
{This paper was written under the support 
of RSF (Grant no~24-44-00099).}
}

\author{Vladimir Kanovei\thanks
{Institute for Information Transmission Problems
(Kharkevich Institute) of Russian Academy
of Sciences (IITP), Moscow, Russia, 
{\tt kanovei@iitp.ru} } 
\and 
Vassily~Lyubetsky\thanks
{Institute for Information Transmission Problems
(Kharkevich Institute) of Russian Academy
of Sciences (IITP), Moscow, Russia, 
{\tt lyubetsk@iitp.ru} } 
}

\date{\today}
\maketitle


\begin{abstract}
We demonstrate that theories 
$\zm$, 
$\zfm$, $\zfcm$ 
(minus means the absence of the Power Set axiom) 
and $\pad$, $\padm$ 
(minus means the absence of the Countable Choice 
schema) are equiconsistent to each other.
The methods used include the interpretation of 
a power-less set theory in $\padm$ via well-founded 
trees, as well as the G\"odel 
constructibility in the said power-less set theory.
\indent

\end{abstract}

\footnotetext{MSC 03E25, 03E35, 03F35 (Primary), 
03E15 (Secondary).}


%
\np

\setlength{\cftbeforesubsecskip}{2.2pt}

{\def\contentsname{%
\normalsize
\bf
Contents
%
%
\vspace*{-0.50ex}
}
\small
\tableofcontents}



%
\np

\parf{Introduction}
\las{int}


The following 
theorem is {\bf the main result} of this paper. 

\vyk{
Theorem \ref{2t} below is the second main result, along 
with the associated technique of consistency proofs by 
class-forcing over the power-less set theory $\zfcm,$ or 
equivalently, the 2nd order Peano arithmetic $\pad$. 
}

\bte
\lam{1t}
Theories\/ $\padm,$  $\pad$, $\zm,$  
$\zfcm,$  $\zfm$ 
are equiconsistent.
\ete

Here $\pad$, resp., $\padm$ 
is 2nd order Peano arithmetic 
with, resp., without the (countable) $\AC$, 
whereas 
$\zm$ is Zermelo set theory without the 
well-orderability axiom $\WA$,
and $\zfcm /\:\zfm$ are  
Zermelo--Fraenkel set theories resp.\ 
with/without $\WA$, and 
all three of them without the {\axf Power Set} 
axiom. 
(See the exact definitions in 
Section~\ref{prePA} 
related to 2nd order Peano arithmetic, 
and in Section~\ref{2prel} related to power-less 
set theories.) 

The theorem has been known since at least 
late 1960s, see \eg\ \cite{Kr}. 
However apparently no 
self-contained and more or less complete 
proof has ever been published. 
(See a brief discussion in {\em Mathoverflow\/} 
around \cite{mof}.) 
The purpose of this paper is to finally present 
such a proof. 

The proof of Theorem~\ref{1t} consists of two 
parts. 
As the {\ubf first part}, we consider 
(Sections~\ref{2prel} and \ref{dev})  
a theory $\zmc$, 
which extends $\zm$ 
by 1) the existence of transitive closures, 
2) the \coun\ axiom, 
and 3) an axiom saying that any \wef\ 
relation on $\om$ admits a transitive model. 
This is a subtheory of $\zfm+\coun$ 
strong enough to prove {\axf Replacement} in 
case the range of the function to exist is a 
transitive class (Lemma~\ref{mz12}). 

Our {\ubf first key result} (Theorem~\ref{it1}) 
provides interpretations of $\zmc$ in $\padm$, 
and of $\zfcm$ in $\pad$,
obtained by using \wef\ subtrees of $\om\lom$ 
as the domain of interpretation. 
This is a well-known method, presented in 
\cite{aptm,marekA,zbi71,zbi81} among others, 
as well as in Sections 3--6 of Chapter VII 
of Simpson~\cite{simp}, 
and in 
\cite{kam} \poo\ 2nd order set theory. 
The tree structure $\bI$, related to this 
interpretation, is defined and studied in 
Section~\ref{dim}. 
The ensuing 
Corollary~\ref{red1} claims the existence 
of two groups of mutually interpretable 
and equiconsistent theories, which include, 
in particular, 
resp., $\padm$, $\zm,$ $\zmc$ 
(group 1)
and $\pad$, $\zfm,$ $\zfcm$ 
(group 2). 

The {\ubf second part} of the proof presents 
an interpretation of $\zfcm$ in $\zmc$,   
contained in the following theorem, 
our {\ubf second key result} here:

\bte
[$\zmc$]
\lam{2red}
The following set or class 
satisfies\/ $\zfcm:$\vim
$$
\Lz\ \ =\ \ \ \left\{
\bay{llr}
\rL\;,& \text{\rm in case $\omil$ does not exist}&
\hspace*{10ex}
{\rm(a)}\\[0.8ex]
\rL_\Om=\bigcup_{\al<\Om}\rL_\al\;,
& \text{\rm in case $\omil=\Om$ exists}&
\hspace*{10ex}
{\rm(b)}
\eay
\right.
$$
\ete

\vyk{
\begin{subequations}
$\Lz$=\left\{
\begin{align}
xx & = 0,
                && \forall\  \underline{z}_h \in\mathbb{Q}^k(K), \\[1mm]
yy,
                && \forall\  p \in \mathbb{P}^k(K),              \\[1mm]
zz  & = 0,
                && \forall\  \mu\in\mathbb{P}^k(e).
\end{align}
\right.
\label{EE}
\end{subequations}

\eqref{EE}
}

\vyk{
\bte
[$\zmc$]
\lam{2red}
The class\/ 
$\Lz=\bigcup_{\al\in\Ord\text{\rm\ is countable in }\rL}
\rL_\al$
satisfies\/ $\zfcm.$ 
We may note here that\vimh
\ben
\raenu
\itlb{2red1A}%
if\/ $\omil$ does not exist then\/ $\Lz=\rL\,;$

\itlb{2red1B}%
if\/ $\omil=\Om$ exists then\/ 
$\Lz=\rL_\Om\,.$
\een
\ete
}

Theorem~\ref{2red} provides an interpretation 
(namely, $\Lz$) of $\zfcm$ in $\zmc$, hence   
connects groups 1 and 2 above, and thereby 
{\ubf implies 
the equiconsistency result of Theorem~\ref{1t}.}  
This interpretation is close to an interpretation 
defined by Simpson \cite[VII.4]{simp}. 
We review some other interpretations, including 
an early one defined in \cite{mz81}, in 
Section~\ref{2om}.

\vyk{
However the theorem has the additional 
advantage of giving a 
\rit{transitive ``standard''} 
(\ie, with the true membership) interpretations   
of $\zfcm$ 
in $\zmc$, a theory 
essentially weaker than $\zfcm.$ 
}

Regarding the class $\rL$ as a whole 
(without a possible reduction to $\rL_\Om$ 
as in Theorem~\ref{2red}(b)), 
we prove the following:

\bte
[$\zmc$]
\label{3t}
\ben
\Renu
\itlb{3tA}%
$\rL\kpo$ 
satisfies\/ $\pad$. 

\itlb{3tB}%
$\rL$ itself satisfies the axiom of\/ 
{\axf Separation}.
\een
\ete

On the other hand, 
$\rL$ does not necessarily satisfy $\zfcm$ 
under $\zmc$ (see Example~\ref{noco}), 
hence option (b) of 
Theorem~\ref{2red} definitely 
cannot be abandoned. 
On the other hand, Theorem~\ref{3t}\ref{3tA} 
asserts that 
$\rL\kpo$ outright satisfies $\pad$ under $\zmc$.

\vyk{
In this corollary, the constructibility is 
meant by referring to a certain $\lpad$  
formula explicitly defined by 
Addison~\cite{add2,add1}. 
}

Theorem~\ref{2red} is proved  
in Sections~\ref{mas}, \ref{fura1}  
on the basis of G\"odel's  
constructibility, developed in 
Sections~\ref{Lopr}, \ref{abs} 
in the context of  $\zmc$.  
%
\vyk{
It is an interesting issue here, 
demonstrated by Example~\ref{noco}, 
that $\zmc$ is not strong enough to 
eliminate case \ref{2red1a}, by proving  
it outright that the whole class $\rL$ 
satisfies $\zfcm$ under $\zmc$. 
}%
Section~\ref{mas} contains Theorem~\ref{z37mz}, 
a key result saying that, in $\zmc$ and  
under certain conditions, 
a class of the form $K=\bigcup_{\al\in\Om}\rL_\al$ 
satisfies $\zfcm,$ leading to the proofs of 
Theorems~\ref{2red} and \ref{1t} 
in Section~\ref{fura1}. 

Then we prove Theorem~\ref{3t}   
in Sections~\ref{fura2}, \ref{fura3}.

The ensuing Corollary~\ref{2redC} states that, 
under $\padm$, 
$\rL\kpo$ satisfies\/ $\pad$. 
Saying it differently, $\rL\kpo$ is an 
interpretation of\/ $\pad$ in\/ $\padm$. 


\vyk{
Theorem~\ref{2red} leaves open a question whether 
$\rL$ itself satisfies $\zfcm$ under the $\zmc$ 
axioms in the full universe. 
Example~\ref{noco} will show that the answer is, 
generally speaking, in the negative.  
Nevertheless
}


\vyk{
\bte
[$\zmc$]
\lam{tL}
$\rL$ itself satisfies the axiom of\/ 
{\axf Separation}.
\ete
}

Our proof of Theorem \ref{1t} leaves open the 
following question: is there a way to interpret 
$\pad$ in $\padm$ avoiding substantial use of set 
theoretic concepts and methods such as 
constructibility? 
A possible approach to this goal, based on the 
ramified analytical hierarchy, is outlined in 
Section~\ref{rah}.

Overall, this is a research and review article, 
the purpose of which is to present 
proofs of such principally important results, as 
Theorems \ref{1t}, \ref{2red}, \ref{3t}, 
in a rather self-contained 
and reader-friendly form.

\vyk{
Analytical representation of G\"odel's constructibility 
is well-known since 1950s, see \eg\ Addison~\cite{add1}, 
Apt and Marek \cite{aptm}, and Simpson's book \cite{simp}. 
This raises {\bf the problem of the consistency} of 
(the analytical form of) $(\ast)$  under the assumption that 
only the consistency of $\pad$ as a premise is available, 
rather that the (much stronger) consistency of $\ZFC$. 
{\gom This is why we consider and solve this problem in our 
paper.}

The working technique of such a transformation of the 
consistency results related to $\ZFC$ to the basis of 
$\pad$ is also rather well-known since some time ago. 
}


\parf{Preliminaries: 2nd order arithmetic} 
\las{prePA}

We recall that  
second order 
Peano arithmetic $\pad$ is a theory in the 
language $\lpad$ with two sorts of variables, 
\imar{lpad}%
\index{LPA2@$\lpad$, language}%
for natural numbers and for sets of them. 
 
We'll use $j,k,m,n$ for variables over $\omega$ and 
$x,y,z$ for variables over $\pws\om$, reserving capital 
letters for subsets of $\pws\om$ and other sets. 

The axioms of $\pad$ are the Peano axioms 
for numbers plus the following:%
\index{theory!PA2@$\pad$}%
\index{PA2@$\pad$, theory}%
\imar{pad}%
%
\bde
\item[\axf Induction:]
$\kaz x\, \big(
0\in x\land\kaz n\,(n\in x\imp n+1\in x) 
\imp \kaz n\,(n\in x)
\big).$ 

\item[\axf Extensionality] for sets: 
$\kaz x,y\,\big(\kaz k\,(k\in x\eqv k\in y) 
\imp x=y \big)$.

\item[\axf Comprehension]  $\CA$:
$\sus x \,\kaz k\,(k\in x\eqv\Phi(k))$ --- 
for every formula $\Phi$ in which $x$ 
does not occur, and in $\Phi$ 
we allow parameters, 
that is, free variables other than $k$.

\item[\axf Countable Choice] 
$\ACw$: 
\imar{ACw}%
$\kaz n\,\sus x\,\Phi(n,x)\imp
\sus x\,\kaz n\,\Phi(n,(x)_n))$ ---  
for every formula $\Phi$ with 
parameters, where  
%
%
$(x)_n=\ens{j}{\asko n j\in x}$, and\break 
$\asko nj=2^n(2j+1)-1$ is a standard bijection 
\imar{asko nj}%
\index{$\asko nj$, standard bijection}%
$\om\ti\om$ onto $\om$.

\ede

The theory $\pad$ is also known as $A_2$ 
(see \eg\ an early survey \cite{aptm}), 
as $Z_2$ (in \cite{schindt} or elsewhere). 
See also \cite{HFuse81,Kr,simp}.
Let $\padm$ 
\index{theory!PA2m@$\padm$}%
\index{PA2-@$\padm$, theory}%
\imar{padm}%
be $\pad$ sans $\ACw$.%
\vom

{\ubf Coding in 2nd order arithmetic}.  
\label{kod}
It can be viewed as a certain disadvantage   
that $\padm$ doesn't 
explicitly treat such objects as pairs, tuples 
and finite sets of numbers, as well as trees of 
tuples at the next level. 
However, these and similar 
(and in fact even more complex) mathematical 
objects can be effectively encoded as single 
natural numbers or sets of them. 
We refer to \cite{simp}, Chap.\,I, and especially 
Section II.2, with respect to many examples. 

Recall that $\Seq=\om\lom$ is the set of all 
tuples (finite sequences) of numbers in $\om$. 
If $s\in\Seq$ and $j<\om$ then $s\we j\in\Seq$ 
is obtained by adjoining $j$ as the rightmost 
term. 
Let $\lh s$ denote the length (the number of terms). 

Let $s_0=\La$ (the empty tuple), and, 
by induction, if $n=\asko mj+1\ge1$ then 
$s_n=s_m\we j$. 
Clearly $\Seq=\ens{s_n}{n<\om}$, and in fact 
$n\mto s_n$ is a bijection onto $\Seq$. 
Subsequently $n=n(s)$ is viewed as  
\rit{the code} of any $s=s_n\in\Seq$,  
and a set $x\sq\om$ is viewed as  
\rit{the code} of   
$\ens{s_n}{n\in x}\sq \Seq$.
Following \cite[esp.\ II.2]{simp}, 
this enables us to freely consider tuples and sets 
of them as if they properly exist, but still on the 
basis of $\padm$.  

Similarly,  still based of $\padm$, 
we can treat sets 
$X\sq\om\ti\om$, $H\sq\Seq\ti\Seq$, 
and the like, as 
properly existing.  

Infinite (and also finite) 
sequences of subsets of $\om$ 
are within reach in $\padm$ as well, because each 
$x\sq\om$ is a code of the infinite sequence of 
sets 
$(x)_n=\ens{j}{\asko n j\in x}$ 
(see the formulation of $\ACw$ above).
And so are \eg\ infinite sequences of subsets 
of $\Seq$.

\parf{Preliminaries: 
intermediate power-less set theory} 
\las{2prel}

{\it The power-less set theory $\zfcm$\/} 
\index{theory!ZFC-@$\zfcm$}%
%
%
is a subtheory of $\zfc$ obtained 
so that:
\ben
\Renu
\itlb{zfcm1}%
the Power Set axiom {\PS} is excluded --- 
the upper minus in $\zfcm$ symbolizes  
the absence of {\PS};

\itlb{zfcm2}%
the usual set-theoretic Axiom of Choice {\AC} 
of $\ZFC$ is removed 
(as it does not work properly without {\PS}), 
and instead the \rit{well-orderability axiom\/} 
{\WA} 
is added, 
which claims that every set can be well-ordered;

\itlb{zfcm3}%
Separation $\Sep$ is preserved, but
the Replacement schema $\Rep$
(too weak in the absence of {\PS}) 
is substituted with the {\em Collection\/} 
schema:\vom 

$\Col:$  
\imar{Col}%
$\kaz X\, 
\big(\kaz x\in X\,\sus y\,\Phi(x,y)\imp 
\sus Y\,\kaz x\in X\,\sus y\in Y\,\Phi(x,y)\big)$. 
\vom

\noi
Note that ${\Col+\Sep}\imp\Rep$. 
\een

See \cite{AG,gitPWS,gitm} for  
a comprehensive account of main features of $\zfcm.$ 

See \cite{jechmill}, \cite{zarRC}, 
\cite[Sect.\,2]{devL} or elsewhere 
for different but equivalent formulations of 
Collection, as \eg\ the following 
form in \cite[Chap.\,6]{jechmill}:\vtm

\quad$\Col'$: 
$\kaz X\,\sus Y\, \kaz x\in X\,
\big(\sus y\,\vpi(x,y)\imp 
\sus y\in Y\,\vpi(x,y)\big)$.\vtm 

\noi
This is apparently stronger than $\Col$ above, 
yet in fact $\Col'$ is a 
consequence of ${\Col}$, for let 
$
\Phi(x,y) := \vpi(x,y) \lor
\big(
y=0\land \neg\,\sus y\,\vpi(x,y)
\big)
$
in \Col.

\bit
\item\msur
$\zfm$ is $\zfcm$ without the 
\index{theory!ZFm@$\zfm$}%
%
\imar{zfm}%
well-orderability axiom $\WA$;\\[0.5ex] 
\vyk{
\item\msur
$\zcm$ is $\zfcm$ without the 
\index{theory!ZCm@$\zcm$}%
\imar{zcm}%
Collection schema $\Col$; 
}
%
$\zm$ is $\zfm$ without the 
\index{theory!Z-@$\zm$}%
\imar{zm}%
Collection schema $\Col$. 
\eit

\bdf
\lam{zmad}
Let $\zmc$ be  
$\zm$ plus the following three 
extra axioms: 
\index{theory!FTMw@$\zmc$, intermediate}%
\index{FTMw@$\zmc$, intermediate theory}%
\index{$\zmc$, intermediate theory}%
\imar{zmc}%
\bde
\vyk{
\item[{\axf\Pzu}, \pclo:]%
\imar{pclo}%
\index{axiom!pairc@\pzu, \pclo}%
for any set $X$ there is 
a {\em\pz\/} superset $Y\qs X,$ 
\index{finite-closed}%
\index{set!finite-closed}%
\ie, it is required that if $x\sq Y$ is 
finite then $x\in Y\,;$
}
%
\item[{\axf Transitive superset}, \trsu:]
\index{axiom!tranc@transitive superset, \trsu}%
\imar{trsu}%
for any $X$ there is a  
transitive superset $Y\qs X.$ 

It follows by {\axf Separation} that 
the transitive closure $\TC X$ 
properly exists. 
Recall that $Y$ is {\em transitive\/} 
\index{transitive}%
\index{transitive closure $\TC X$}%
\index{$\TC X$}%
\index{set!transitive}%
if $\kaz x\,\kaz y\,(x\in y\in Y\imp x\in Y)$, 
and the {\em transitive closure\/} is the 
intersection of all transitive supersets.  
%

\item[{\axf Mostowski Collapse, \moco}:]%
\index{axiom!most@Mostowski collapse, \moco}%
\imar{moco}%
\sloppy
any \wf\ relation $A\sq\om\ti\om$ 
admits a transitive set $X$ and 
a function 
$\xet:\obl A=\dom A\cup\ran A\text{\rm\ onto }X$, 
satisfying   
(*)  $\xet(k)=\ens{\xet(j)}{j\mbA k)}$, for all 
$k\in \obl A$.
\index{$\xet_A$}%
\imar{xet A}%
By standard arguments, the map $\xet=\xet_A$ and 
the set $X=\liz A$ are 
\imar{liz A}%
\index{$\liz A$}%
unique.

A binary relation $A\sq\om\ti\om$ 
is {\em \wf\/}, or {\em regular\/}, iff any set
$\pu\ne X\sq\obl A$ contains an element $y\in X$ 
such that $\kaz x\in X\,\neg\,(x\mbA y)$.


\item[{\axf \coun}:]%
\index{axiom!count@\coun}%
\imar{coun}%
$\kaz x\,\sus f\,(f:x\to\om\text{ is 1--1})$, 
\ie, all sets are at most countable$.$
\ede
The name $\zmc$ reflects the initial  
letters of the additional axioms.
\edf

The axiom \moco\ was introduced in \cite{mz81}; 
it is called Axiom Beta in 
\cite[Def.~VII.3.8]{simp}. 
See \cite[Theorem 6.15]{jechmill} for a proof 
of \moco\ in $\ZF$. 

Applying \moco\ for $A={{\in}\res D}$, we immediately 
get:

\bcor
[$\zmc$, transitive collapse]
\lam{c12}
Let\/ $D$ be any set. 
There is a unique transitive set\/ $X$ and a 
unique collapse map\/ $\ta:D\text{\rm\ onto } X$ 
satisfying\/ 
$\ta(x)=\ens{\ta(y)}{y\in x\cap D}$ for all \/ 
$x\in D$.\qed
\ecor

\vyk{
\bre
[$\zm$]
\lam{Dtran}
$(1)$ 


\vyk{
$(2)$ 
If a set $X$ is {\em \pz\/} as in \pclo\ 
then $X$ also is closed under ordered pairs 
$\ang{x,y}$ and generally tuples 
$\ang{x_1,\dots,x_n}$ of any length. 
}

$(3)$ 
The {\em transitive closure\/} $\TC X,$ 
resp., the {\em finite closure\/} $\PC X,$ 
of $X$  
\index{set!transitive closure, $\TC\cdot$}%
\index{transitive closure, $\TC \cdot$}%
\index{set!finite closure, $\PC \cdot$}%
\index{finite closure, $\PC \cdot$}%
is the \dd\sq least transitive, 
resp., \pz\ superset $Y\supseteq X$ (if exists).
\ere

\bre
[$\zmc$]
\lam{nodef}
Dealing with digraphs $A$ in $\wfe$ \etc\
in set-theoretic environment, as in \moco,
it is {\ubf not} assumed  
by default that $A\sq\om\ti\om$. 
We'll explicitly indicate this assumption 
whenever need be.
\ere

}

\ble
[$\zmc$]
\lam{susA}
For any set\/ $X$ there is a \wf\ relation\/ 
$A\sq\om\ti\om$, such that\/ 
$X=\liz A$. 
\ele

\bpf
Let $T=\TC{X}$, the least transitive set 
with $X\sq T$. 
Let $f:T\to\om$ be an injection, by \coun, 
and $A=\ens{\ang{f(x),f(y)}}{x\in y\in X}$. 
Then $A$ is a \wf\ relation on $\om$, 
$\xet_A=f\obr$ satisfies (*) of 
Definition~\ref{zmad}, 
and $X=\liz A$. 
\epf

\bpro
[not used below] 
\lam{2zma}
$\zmc
\sq \zfm+\coun$.\qed 
\epro

Simpson \cite[VII.3.3 and VII.3.8]{simp}  
considers a related theory 
$\atros$ in the $\in$-language, containing the 
\imar{atros}%
\index{$\atros$, theory}%
\index{theory!$\atros$}%
following axioms: 
\ben
\aenu
\itlb{aeq}%
Axiom of Equality: 
$=$ is an equivalence relation and   
$\in$ is $=$-invariant; 

\itlb{aei}%
Axioms of Extensionality and Infinity 
in their usual forms; 

\itlb{rud}%
Axiom of Rudimentary Closure, which asserts,
for all $u, v, w$, the proper existence of 
$\ans{u, v}$, $u\bez v$, $u \ti v$, $\bigcup u$, 
${\in}\res u$, and the following:\vim 
$$
\bay{rcl}
u\obr &=& 
\ens{\ang{x,y}}{\ang{y,x}\in u}\,,\\[0.5ex]
&& 
\ens{\ang{y,\ang{x,z}}}
{\ang{y,x}\in w\land z\in u}\,,\\[0.5ex]
&& 
\ens{\ang{y,\ang{z,x}}}
{\ang{y,x}\in w\land z\in u}\,,\\[0.5ex]
&& 
\ens{v}
{\sus x\,(x\in u\land v=\imb w{\ans x}}\,.\vim
\eay
$$
\itlb{ar}%
Axiom of Regularity 
in its usual form; 

\itlb{a3}%
Axioms \trsu, \moco, \coun, as in 
Definition \ref{zmad}.
\een

\bre
\lam{AE}
$\zmc\bez{\axf Separation}\sq\atros\sq\zmc$. 
Indeed, regarding the second $\sq$, all operations 
listed in \ref{rud} are properly convergent within 
any transitive finite-subset-closed set. 
Now refer to Lemma~\ref{pc} below.
\ere

\vyk{
\bpf 
Standard proofs of \ref{zma0}, \ref{zmaB}, \ref{zmaC} 
of Definition \ref{zmad}  
can be maintained avoiding the 
{\ubf PS} axiom. 
For instance, given a set $X=X_0$, 
we define a 
sequence of sets $X_n$ by unduction, 
setting $X_{n+1}=\bigcup X_n$, and 
then $Y=\bigcup_nX_n$ is a set by the Replacement 
and Union axioms, 
whereas $Y=\TC X$ holds by construction. 

It takes a bit more effort to establish \ref{zmaC}. 
Let\/ $A\in\wfev$. 
Recall that 
$\elm Ak=\ens{j\in\obl A}{j\mbA k}$.
Say that 
$f$ is an $A$-\rit{suitable map} iff 
$D=\dom f\sq \obl A$ is 
$A$-transitive, and we have 
$f(k)=\ens{f(j)}{j\in\elm Ak}$ 
for each $k\in D$.  
We claim that:
\ben
\fenu
\itlb{2z1}%
any two $A$-suitable maps coincide on the common domain;

\itlb{2z2}%
if $k\in\obl A$ then there exists an $A$-suitable map $f$ 
with $k\in\dom f$.
\een

To prove \ref{2z1}, suppose that $f,g$ are suitable and 
$k\in\dom f\cap\dom g$, and hence 
$\elm Ak\sq \dom f\cap\dom g$. 
Assuming that $f(j)=g(j)$ for all $j\in\elm Ak$, 
we get
\pagebreak[0] 
\blin
$$
f(k)=\ens{f(j)}{j\in\elm Ak}=\ens{g(j)}{j\in\elm Ak}=g(k),
$$
\elin
so that still $f(k)=g(k)$. 
It follows by Lemma \ref{wfin} that $f(k)=g(k)$ 
for all $k$. 
This ends the proof of \ref{2z1}. 
The  proof of \ref{2z2} is similar. 

The union  $F$ of all $A$-suitable maps  
is $A$-suitable by \ref{2z1}, 
whereas $\dom F=\obl A$ by \ref{2z2}. 
We claim that $F$ is 1--1.

Indeed let 
$W=\ens{k\in\obl A}{\kaz n\,(F(n)=F(k)\imp n=k)}$. 
We have to prove $W=\obl A$. 
To make use of Lemma \ref{wfin}, 
assume that $k\in\obl A$ and $\elm Ak\sq W.$ 
Suppose to the contrary that $k\nin W.$ 
Let this be witnessed by $n\in\obl A$, so that 
$n\ne k$ but $F(n)=F(k)$. 
Then  
\bce
$
\ens{F(j)}{j\in\elm Ak}=
\ens{F(m)}{m\in\elm An}
$
\ece
by definition. 
As $\elm Ak\sq W,$ we conclude that 
$
\elm Ak=\elm An 
$, 
and hence $k=n$ as $A$ is extensional. 
The contradiction obtained proves $W=\obl A$. 
Thus $F$ is 1--1. 
It easily follows that ${j\mbA k}\eqv {F(j)\in F(k)}$. 

We finally claim that $\ran F=\tca x$, where 
$x=F(M)$ and $M=\verx E\in\dom F$. 
Indeed if $y=F(j)$, $j\in\obl A$, $j\ne M$, then  
there is a finite $A$-chain 
${j=j_0}\mbA j_1\ldots j_{K-1}\mbA j_K=M$ 
in $\obl A$, hence 
$y=F(j_0)\in F(j_1)\in\ldots\in F(j_{K-1})\in F(M)=x$
and $y\in\tca x$. 
Conversely let $y\in\tca x$,  
hence $y=y_0\in y_1\in\ldots\in y_{K-1}\in x$ 
for some $y_n$, $n<K$. 
Then by reverse induction $y_n=F(j_n)$, 
where $j_n\in \obl A$ and 
${j_0}\mbA j_1\mbA\ldots \mbA j_{K-1}\mbA M$. 
In particular $y=F(j_0)\in\ran F$. 
Thus $\ran F=\tca x$.
We conclude that $F$ witnesses \ref{zmaC} with 
$x=F(M)$ and $M=\verx E\in\dom F$. 
\epf
}

\vyk{

\parf{Isomorphism and rigidity of trees}
\las{isr}

Let $A,B$ be two trees. 
We write $A\cong B$ iff there is 
\index{$A\cong B$}%
\index{tree!isomorphism, $A\cong B$}%
an $A,B$-isomorphism, \ie, a bijection 
$H:\obl A\onto\obl B$, satisfying 
${u\mbA v}\iff{H(u)\mbB H(v)}$. 
\pagebreak[0]


{\ubf Important note.}
We may observe that the digraph-extended 2nd-order
arithmetic $\padmx$ 
allows essentially the same tools of dealing with 
digraphs on $\om$ (\ie, subsets of $\om\ti\om$) 
as $\zmx$ (with Cartesian products) does. 
In particular, 
}

\vyk{
\bdf
\lam{HAB}
Given two digraphs $A,B$, define 
a map $H_{AB}$ as follows. 
\index{digraph!map $H_{AB}$}%
If $u\in\obl A$, $v\in\obl B$, and 
$\con Au\cong\con Bv$ then put $H_{AB}(u)=v$. 
(This is well-defined since the uniqueness of $v$ 
holds by Lemma~\ref{rig}.) 
\edf

\vyk{
Put 
\blin
\bul
{DD}
{
D=\ens{u\in\obl A}
{\sus v\in\obl B\,(\con Au\cong\con Bv)}\,.
}
\elin
If $u\in D$ then let $H_{AB}(u)$ be the unique 
(by Lemma~\ref{rig}) 
element $v\in \obl B$ satisfying 
$\con Au\cong\con Bv$. 
}

\bcor
[$\padmx$, of Lemma~\ref{rig}]
\lam{hab}
Assume that\/ $\pu\ne A,B\in\wfe$. 
Then\/ $X=\dom H_{AB}\sq\obl A$ 
is\/ $A$-transitive, 
$Y=\ran H_{AB}\sq\obl B$ is\/ $B$-transitive, 
and\/ $H_{AB}$ is a bijection\/ $X\onto Y$ which 
witnesses\/ ${A\res X}\cong {B\res Y}$, 
so that\/ ${u\mbA v}\iff{H_{AB}(u)\mbB H_{AB}(v)}$, \ 
for all\/ $u,v\in X\;.$ \qed
\ecor 
}

\vyk{
\itlb{hab0}\msur%
the sets\/ $X,Y$ contain resp.\ 
the\/ $A$-minimal element\/ $a=\min A$\/ 
{\em(see Remark~\ref{mine})}  
and the $B$-minimal element $b=\min B,$  
and $H_{AB}(a)=b;$   

\itlb{hab2}\msur%
$H_{AB}$ is a bijection\/ $X\onto Y$ which 
witnesses\/ ${A\res X}\cong {B\res Y}$, 
so that\/ ${u\mbA v}\iff{H_{AB}(u)\mbB H_{AB}(v)}$, \ 
for all\/ $u,v\in X\;.$ 
  
\itlb{hab3}%
if\/ $u\in\obl A\bez X$ then no\/ $v\in\obl B$ 
satisfies\/ $\con Au\cong\con Bv$, \break
if\/ $v\in\obl B\bez Y$ 
then no\/ $u\in\obl A$ 
satisfies\/ $\con Au\cong\con Bv\,.$ 
}


\parf{Development of the 
intermediate power-less theory} 
\las{dev}

We proceed with a few simple results in $\zmc$ 
hardly available in $\zm.$ 

Let a \rit{class-function} be a (definable) class 
\index{class-function}%
which satisfies the standard definition of a 
function 
(\ie, consists of sets that are ordered pairs, \etc).

\vyk{
\ble
[$\zmc$]
\lam{tc}
For any set\/ $X,$ the transitive closure\/ $\TC X$ 
properly exists as a set.
\ele

\bpf
Note that  
$\TC X=\ens{y\in Y}{\tau(y)}$, 
where $Y\qs X$ is a transitive superset 
(exists by \trsu) and  
$\tau(y)$ says: {\em there is a finite 
sequence 
$y=y_n\in y_{n-1}\in\ldots\in y_1\in y_0\in X$ 
of sets $y_i\in Y.$} 
(Including the case $n=0$ and $y\in X.$) 
\epf
}

\ble
[$\zmc$]
\lam{mz12}
Let\/ $F$ be a class-function, $D=\dom F$ any set. 
Then\/ $F$ and the image\/ 
$R=\imd FD=\ens{F(x)}{x\in D}$ are sets in each of the 
two cases$:$ $(1)$ $R$ is transitive, $(2)$ 
there is a set\/ $Y$ such that\/ $R\sq \pws Y.$  
\ele

\bpf 
(1)
By \coun\ we can \noo\ assume that $D\sq\om$. 
We can also assume that $F$ is 1--1, 
for otherwise replace $D$ by the set 
\bce
$D'=\ens{k\in D}{\kaz j\in D\,(j<k\imp F(j)\ne F(k))}$.
\ece
Then the relation 
$A=\ens{\ang{j,k}}{j,k\in D\land F(j)\in F(k)}$ 
is \wf\ as isomorphic to ${\in}\res R$. 
On the other hand, by \moco, $A$ is isomorphic 
to ${\in}\res Y$ where $Y$ is a transitive \rit{set}. 
It follows that $Y$ and $R$ are $\in$-isomorphic, 
and hence $R=Y$ is a set. 
Finally $F\sq X\ti R$ is a set by Separation. 

(2) 
We \noo\ assume that $Y$ is transitive, by 
\trsu. 
We can assume as well that $D\cap Y=\pu$, 
otherwise put $D'=D\ti\ans Y$ and change $F$ 
accordingly. 
Under these assumptions, put $D_1=D\cup Y$ and 
extend $F$ to $F_1$ by the identity on $Y.$ 
Then the image $\imd{F_1}{D_1}=R\cup Y$ is transitive,  
hence, a set by (1). 
Now $R\sq \imd{F_1}{D_1}$ is a set by \Sep.
\epf

A set $Y$ is called {\em finite-subset closed} 
\index{finite-subset closed}%
\index{set!finite-subset closed}%
if $\kaz z\sq Y\,(z\text{ finite }\imp z\in Y).$
For any set $X$, let the {\em\pzu\/} $\PC X$ be 
\index{$\PC X$}%
\index{finite closure}%
the least finite-subset closed superset $Y\qs X$, 
if exists.

\ble
[$\zmc$]
\lam{pc}
For any set\/ $X,$ $\PC X$ 
properly exists.
\ele

\bpf
To handle the case $X=\om$, let $p_k$ be $k$th prime, 
so $p_1=2$, $p_2=3$, and so on. 
Let 
$A=\ens{\ang{k,n}}{k\ge1\land p_k\text{ divides }n}$.
Then $\obl A=\om\bez\ans0$, $A$ is \wf\ 
(since $k\mbA n\imp k<n$), and (\dag) for any finite 
$u\sq\obl A$ there is $n\in\obl A$ satisfying 
$u=\ens{k}{k\mbA n}$. 
By \moco\ there is a map 
$\xet:\obl A
\text{\rm\ onto a transitive set }R$, 
satisfying   
(*)  $\xet(n)=\ens{\xet(k)}{k\mbA n)}$, for all 
$n\in \obl A$.
Then easily $R=\PC\om$ by (\dag).

To handle the general case, we may assume that 
$X$ is transitive, by \trsu.
Let $h:\om\text{ onto }X$, by \coun. 
Then $h$ can be extended to a class-map $H$ 
defined on the bigger set 
$R=\PC\om$ so that $H\res\om=h$ and if 
$u\in R\bez\om$ then $H(u)=\ens{H(n)}{n\in u}$.
Then $\ran H=\PC X$ (so far a class), 
and hence $\ran H$ is transitive as so is $X$. 
It follws by Lemma~\ref{mz12} that both $H$ and 
$\ran H=\PC X$ are proper sets.
\epf

\ble
[$\zmc$]
\lam{capr}
Let\/ $U,V$ be any sets. 
Then\/ $U\ti V$, 
$\pfin U$, $U\lom$ properly exist (as sets).
\ele
\bpf 
$X=U\cup V=\bigcup\ans{U,V}$ 
is a set by $\zm.$ 
Now $\PC X$ is a set by Lemma~\ref{pc}, 
hence $U\ti V\sq \PC X$ is a set by 
\Sep. 
To prove the other claims, 
note that 
$\pfin U\yi U\lom\sq \PC U$ and use 
Lemma~\ref{pc} and \Sep.  
\epf

Thus $\zmc$ proves the existence of Cartesian 
products. 
Note that $\zm$ does not prove even the existence 
of $\om\ti\om$!  

\ble
[$\zmc$]
\lam{ord}
Let\/ $E$ be a strict \weo\ of a set\/ $U.$  
Then there is an ordinal\/ $\la$ and an 
order isomorphism of\/ $\stk UE$ onto\/ 
$\stk\la\in$.
\ele
\bpf 
By \coun\ we can \noo\ assume that $U\sq\om$. 
Then $E$ is a \wf\ relation with $\obl E\sq\om$. 
Apply $\moco$. 
Then $\la=X$ is a transitive set well-ordered 
by $\in$, that is, an ordinal.
\epf

\bcor
[$\zmc$]
\lam{ordC}
If\/ $\al,\ba$ are ordinals then there exist 
(as sets)   
ordinals\/ $\al+\ba$, $\al\cdot\ba$, $\al^\ba$. 
{\em(In the sense of the ordinal arithmetic.)}
\ecor
\bpf 
We have to define well-ordered sets which 
represent the mentioned orders. 
For instance, the Cartesian product $\al\ti\ba$ 
(a set by Lemma~\ref{capr}), 
ordered lexicographically, represents $\al\cdot\ba$. 
The exponent $\al^\ba$ is represented by 
the set 
$W=\ens{f:D\to\al\bez\ans0}
{D\sq\ba\text{ is finite}}$, 
ordered lexicographically, with the 
understanding that each $f\in D$ is by default 
extended by $f(\xi)=0$ for all $\xi\in\ba\bez D$. 
Note that $W\sq\PC{\ba\ti\al}$ is a set by 
Lemma~\ref{pc}. 
\epf

\parf{The set theoretic tree hull 
over second-order PA}
\las{dim}

Following \cite[VII.3]{simp}, we consider the 
collection $\wft$ of all \wf\ trees 
$\pu\ne T\sq\Seq=\om\lom$. 
Recall that 
\bit
\item\msur
$\La$ is the empty tuple, 
$\ang k$ is the tuple with $k$ as the single term;  
\index{$\La$, the empty tuple}%

\item\msur
$T\sq\Seq$ is a tree iff 
$s\we j\in T\imp s\in T$;  

\item\msur
$T$ is \wf\ iff 
$\neg\,\sus g:\om\to\om\,\kaz m\,(g\res m\in T)$;

\item\msur
$s\we j$ is obtained by adding $j\in\om$ to $s\in\Seq$ 
as the rightmost term, and if $s,t\in\Seq$ then 
$s\we t\in\Seq$ is the {\em concatenation}; 
\index{concatenation, $s\we t$}%
\index{$s\we t$}%
\index{$s\we j$}%

\item%
if $T$ is a tree and $s\in T$ then put  
$T^s=\ens{t\in\Seq}{s\we t\in T}$; thus $T^s$ is a tree 
\index{$T^s$}%
as well, and if $T$ is \wf\ then so is $T^s$.
\eit

\bdf
[$\padm$]   
\lam{bis}
Let $S,T\in\wft.$ 

A set $H\sq S\ti T$ is 
\index{tree!bisimulation}%
\index{tree!$A\cong B$}%
an $S,T$-\rit{bisimulation}, iff, for all 
$s\in S$ and $t\in T$, 
\blin
\bul
{bis1}
{{s\mbH t}\,\iff\, 
\bay[t]{l}
\phantom{\land\;\;}
\kaz s'=s\we j\in S\;\sus t'=t\we k\in T\;(s'\mbH t')
\;\;\land\,\\[0.5ex]
\land\;\;
\kaz t'=t\we k\in T\;\sus s'=s\we j\in S\;(s'\mbH t')\,.
\eay}
\elin

Define $S\cong T$ iff there is 
\index{$A\cong B$}%
an $S,T$-{bisimulation} $H$ such that $\La\mbH\La$.


Define $S\inh T$ iff 
\index{$A\inh B$}%
\index{tree!$A\inh B$}%
\index{relationinh@relation $\inh$}%
$S\rah T^u$ for some $u\in T$ with $\lh u=1$. 

{\em The structure\/} 
$\bI={\stk{\wft}{\rah,\inh}}$ 
\index
{structureV@structure $\bI=\stk{\wft}{\rah,\inh}$}%
is considered in $\padm$. 
\index{relationinh@relation $\inh$}%
\index{relationrah@relation $\rah$}%

The $\bI$-{\em interpretation\/} $\jnt\Phi$ 
of an\/ $\in$-formula $\Phi$ 
(with parameters in $\wft$)  
is naturally 
\index{PhiS@$\jnt\Phi$, interpretation}%
\index{$\jnt\Phi$, interpretation}%
defined in the sense of 
interpreting\/ $=,{\in}$ 
as resp.\ ${\rah},{\inh}$, and relativizing 
the quantifiers to $\wft$. 
Thus \eg\ $\jnt{x=y}$ is $x\rah y$.
\edf

Note that the bisimulation relation $\cong$ between 
trees in $\wft$, and subsequently the derived 
relation $\inh$ as well, are naturally formalized 
in $\padm$ in the frameworks of the approach 
based on coding see Section~\ref{kod}. 
It follows that, for any $\in$-formula $\Phi$ 
with parameters in $\wft$, the 
$\bI$-{interpretation\/} $\jnt\Phi$ 
of is a $\lpad$-formula.  

\bte
[$\padm/\pad$]   
\lam{it1}
$\bI$ is a well-defined structure$:$ 
$\rah$ is an equivalence on\/ $\wft,$  
$\inh$ is a binary relation on\/ 
$\wft$ invariant \poo\/ $\rah$. 

Moreover\/ $\bI$ 
satisfies resp.\ $\zmc/\zfcm.$
In other words, if\/ $\Phi$ is an axiom of\/ 
$\zmc$, resp., $\zfcm$ then\/ $\jnt{\Phi}$ 
is a theorem of resp.\/ $\padm$, $\pad$. 
%
\ete

This theorem is a version of the interpretation 
results known since at least Kreisel~\cite{Kr} 
and published somewhat later  
in \cite{aptm,marekA,zbi71}    
or elsewhere. 
In particular, the $\pad$ part of the theorem 
is essentially Theorem 5.5 in \cite{aptm}. 
The $\padm$ part is close to 
Theorem 1.1 and Corollary 1.1 in 
\cite{marekA}. 

\bpf
Besides the papers cited above, the bulk of the 
theorem was established in \cite[VII.3]{simp}. 
Namely, using just $\text{\rm\ubf ATR}^0$ as the basis 
theory (which is a small part of $\padm$), 
Lemma~VII.3.20 in \cite{simp} proves that if 
$\Phi$ is an axiom of $\atros$ 
then $\jnt{\Phi}$ is a theorem of 
$\text{\rm\ubf ATR}^0$ 
(and then of $\padm$ as well). 
Thus, to prove the $\padm$ part of Theorem~\ref{it1},  
it suffices to check \Sep\ in~$\bI.$  

{\ubf Arguing in $\padm,$} assume that $S\in\wft$, 
$X=\ens{k}{\ang k\in S}$, and $\Phi(x)$ is an 
$\in$-formula with parameters in $\wft$ and 
with $x$ as the only free variable. 
Trees of the form 
$S^k=\ens{t\in\Seq}{k\we t\in S}$, $k\in X$, 
belong to $\bI$ and are the only (modulo $\rah$)
$\inh$-elements of $S$ in $\bI$. 
Now, using the $\padm$ {\axf Comprehension}, 
we let $Y=\ens{k\in X}{\jnt{\Phi(S^k)}}$. 
The set $T=\ans\La\cup \ens{t\in S}{t(0)\in X}$ 
is a tree in $\wft$.  
We claim that $\jnt{T=\ens{x\in S}{\Phi(x)}}$. 

Indeed assume that $C\in\wft$, $C\inh S$, and 
$\jnt{\Phi(C)}$. 
Then $C\rah S^k$ for some $k\in X$, 
so that $\jnt{\Phi(S^k)}$ holds, 
and hence $k\in Y.$ 
It follows that $C\rah T^k=S^k\inh T.$  
The proof of the inverse implication is similar.

Finally prove {\ubf the $\pad$ part} of the theorem. 
{\ubf Arguing in $\pad,$} we have to check 
{\axf \Col} in $\bI$. 
Thus let $S\in\wft$ and let $\Phi(x,y)$ be an 
$\in$-formula with parameters in $\wft$, 
satisfying $\jnt{\kaz x\in S\,\sus y\,\Phi(x,y)}$,  
that is, 
\blin
\bul
{kol}
{ 
\kaz A\in\wft\,\sus B\in\wft\,
\big(A\inh S\imp \jnt{\Phi(A,B)}
\big).
}
\elin
But $\inh$-elements of $S$ are, modulo $\rah$, 
trees $S^k=\ens{s\in S}{k\we s\in T}$, 
where $k\in K=\ens{k\in\om}{\ang k\in T}$, and 
only them. 
Thus \eqref{kol} implies: 
\blin
$$
\kaz k\in K\,\sus B\in\wft\,
\big( \jnt{\Phi(S^k,B)}
\big).
$$
\elin
Using $\ACw$ of $\pad$, we get a 
(coded, see Section~\ref{kod}) 
sequence of trees $B_k\in\wft$ with 
$\jnt{\Phi(S^k,B_k)}$ for all $k$. 
Now $T=\ang\La\cup\,\bigcup_{k\in K}k\we B_k\in\wft$, 
and each $B_k$ is an $\inh$-element of $T$. 
Thus we have 
\blin
$$
\kaz k\in K\,\sus B\inh T\,
\big(\jnt{\Phi(S^k,B)}
\big), \text{ \ that is, \ } 
\jnt{\kaz x\in S\,\sus y\in T\,
\Phi(S^k,B)},
$$
\elin
as required. 
\epf

\bcor
[of Theorem~\ref{it1}]
\lam{red1}
Theories\/ $\padm$, $\zm,$ $\zmc$ are 
mutually interpretable and hence equiconsistent. 
Theories\/ $\pad$, $\zfm,$ $\zfcm$ 
are mutually interpretable and 
equiconsistent as well.\qed 
\ecor

Corollary~\ref{red1} is the first part 
of the proof of Theorem~\ref{1t}. 
The remainder of the proof involves the ideas 
and technique of G\"odel's constructibility, 
and {\ubf the goal will be Theorem~\ref{2red}} 
containing an interpretation of 
$\zfcm$ in $\zmc$.

\vyk{

\parf{The definability jump} 
\las{dj}

We argue in $\zmc$ in this Section. 

Let $X$ be a transitive set. 
We are going to legitimize 
the set $\Def X$ of all sets $Y\sq X$, 
$\in$-definable over $X$ 
(\ie, over the structure $\stk X\in$), 
with parameters $p\in X$ allowed --- a basis 
of the standard definition of the constructible 
hierarchy. 
The operation $\Def$ is a well-known $\ZF$ notion,  
see \cite[pp.~157,176]{jechmill}.  
Yet $\zmc$ does not contain some key $\ZF$ 
tools, namely, the {\axf Power Set} axiom, and 
{\axf Collection/Replacement}, except for a small 
piece of \Rep\ by Lemma~\ref{mz12}. 
Hence we have to be careful working in $\zmc$. 

Let $\form$ be the collection of all 
\imar{form}%
\index{$\form$}%
parameter-free $\in$-formulas. 
After a standard g\"odelization, 
any $\vpi\in\form$ is a finite string of 
symbols coded by integers, hence formally 
$\form\sq\om\lom.$ 
But $\om\lom$ is a set in $\zmc$ by Lemma~\ref{pc}. 
We conclude that $\form$ is a recursive subset of  
$\om\lom$ by \Sep. 

Giver a transitive set $X$, we 
let $\form[X]$ be the collection of all 
\imar{form[X]}%
\index{$\form[X]$}%
formulas in $\form$, in which some 
(including all or none) 
free variables is substituted with sets 
$p\in X$ (parameters). 
$\form[X]\sq (X\cup\om)\lom$ is a set 
by \Sep. 

Put 
$\form_\vpi=\ens{\psi}
{\psi\text{ is a subformula of }\vpi}$ --- 
\imar{form vpi}%
\index{$\form_\vpi$}%
this is a finite set.  
Sets $\form_\vpi[X]$ (parameters in $X$ allowed) 
are defined accordingly. 
\imar{form vpi[X]}%
\index{$\form_\vpi[X]$}%
A \rit{Tarski truth set} (TTS) for 
\index{TTS, Tarski truth set}%
$\vpi\in\form[X]$ over a transitive set $X$ is 
any set 
$\tau\sq\ens{\psi\in\form_\vpi[X]}
{\psi\text{ is closed}}$ 
satisfying the standard Tarski conditions:
\ben
\fenu
\itlb{12*}\msur%
$(x\in y)\in\ta$ iff $x\in y$,\\[0.3ex]  
$(\neg\psi)\in\ta$ \ iff \  
$\psi\nin\tau$ --- for closed formulas 
$\psi\in\form_\vpi[X]$,\\[0.3ex]  
$(\psi\land\chi)\in\ta$ \ iff \  
$\psi\in\tau$ and $\chi\in\tau$,\\[0.3ex]  
$(\sus v\,\psi(v))\in\ta$ \ iff \ 
$\psi(x)\in\ta$ for some $x\in X$,
\een
assuming that we have only 
$\neg,\land,\sus$ as logical connectives.

Let $\tar \ta X\vpi$ be the $\in$-formula saying: 
\index{$\tar \ta X\vpi$, formula}%
\bqu\bce
{$X$ is a transitive set, $\vpi\in\form[X]$ is 
a closed formula, 
and $\ta$ is a TTS for $\vpi$ over $X$}.
\ece\equ

\ble
[$\zmc$]
\lam{ttf}
Assume that\/ $X$ is a transitive set, and\/ 
$\vpi\in\form[X]$ is a closed formula. 
Then 
\blin
\bul
{ttfE}
{
\sus \ta \big(\tar \ta X\vpi\land \vpi\in\ta\big)
\iff
\kaz \ta \big(\tar\ta X\vpi\imp \vpi\in\ta\big).
}
\elin
\ele

\bpf
It suffices to establish the existence and 
uniqueness 
of a TTS for $\vpi$ over $X,$ which goes on by a 
routine induction on the length of $\vpi$.
\epf 

\bdf
\lam{dmo}
Define $X\mo\vpi$ iff 
\index{$X\mo\vpi$, definition}%
\index{$\mo$, truth relation}%
$X$ is a transitive set,   
$\vpi\in\form[X]$ is a closed formula, and 
{\em the equivalent formulas in\/ \eqref{ttfE} hold}.
This is the formal truth definition in $\zmc$, 
and both $\vpi$ 
(as a finite string of symbols and 
parameters $p\in X$) 
and $X$ 
are set variables in 
the formula $X\mo\vpi$. 
\edf

\bpro
\lam{ttfC}
Let\/ $\vpi(v_1,\dots,v_n)$ be a 
metamathematical\/ $\in$-formula. 
Then the following claim 
is a theorem of\/ $\zmc:$ 
if\/ $X$ is a transitive set and\/ 
$x_1,\dots,x_n\in X$, 
then\/ 
$\vpi^X(x_1,\dots,x_n)\eqv
\big(X\mo\vpi(x_1,\dots,x_n)\big)$.
\qed
\epro

As usual, here $\vpi^X$ is the formal 
\index{$\vpi^X$, relativization}%
\index{relativization, $\vpi^X$}%
relativization to $X,$  
\cite[p.\,161]{jechmill}.

Given a transitive set $X,$ let 
$\forma[X]=$ all formulas $\vpi(v)\in\form[X]$ 
with exactly one free variable, and let 
$\Def X$ contain all sets of the form 
\index{$\Def X$}%
\blin
\bul
{defX}
{X_{\vpi(v)}=\ens{x\in X}{X\mo\vpi(x)}, \ 
\text{ where } \ 
\vpi(v)\in\forma[X]\,.
}
\elin 
}

\parf{Constructible sets
in the intermediate theory} 
\las{Lopr}

We'll make use of some keynote definitions and 
results related to constructible sets as 
given in \cite[Sect.~VII.4]{simp}. 
We present these results based on $\zmc$,
whereas Simpson works in $\atros$ and in some other
sub-theories of $\zmc$ in \cite{simp}, 
which is not our intention here. 

\ble
[$\zmc$, VII.4.1 in \cite{simp}]
\lam{si1}
Let\/ $X$ be a nonempty transitive set. 
There exists a unique set\/ $\Def X $ consisting 
of all sets\/ $Y \sq X$, definable over\/ 
$X$ by an\/ $\in$-formula with
parameters from\/ $X$.

This set\/ $\Def X$ is obviously transitive, and 
$X\cup\ans X\sq\Def X$. 
%
\qed
\ele

\ble
[$\zmc$, \cite{simp}, Lemma VII.4.2]
\lam{legi}
Let\/ $u$ be a transitive set and\/ $\ba\in\Ord.$  
There is a unique function\/ $f=\jf^u_{\ba}$ such 
\imar{jf}%
\index{$\jf^u_{\ba}$}%
\index{constructing function, $\jf^u_{\ba}$}%
that\/ $\dom f=\ba,$ $f(0)=u,$ 
$f(\al+1)=\Def{f(\al)}$ whenever\/ $\al+1<\ba$, 
and\/ $f(\la)=\bigcup_{\al<\la}f(\al)$ for all 
limit\/ $\la<\ba$. 
\qed 
\ele

The lemma enables us to define 
$\rL_\al[u]=\jf^u_{\al+1}(\al)$ in $\zmc$, 
legitimizing 
the standard definition of relative 
constructible hierarchy
for any set $u\sq\om$:
\pagebreak[0]
\blin
\bul
{LLdef}
{\left.
\bay{lcl}
\rL_0[u] &=& \om\cup\ans u - \text{to keep it transitive},\\[0.7ex]
%
%
\rL_{\al+1}[u] &=& 
\Def{\,\rL_\al[u]}
\text{ \ for all \ }\al\,,\\[0.7ex]
\rL_\la[u] &=& 
\bigcup_{\al<\la}\rL_\al[u] \ \ 
\text{for all limit }\la\,,\\[0.7ex]%
\rL[u] &=& 
\bigcup_{\al\in\Ord}\rL_\al[u] \,= \,
\text{all sets constructible in $u$},\\[0.7ex]%
\rL_\al &=& 
\rL_\al[\pu],\\[0.7ex]%
\rL &=& 
\rL[\pu].
\index{set!constructible, $\rL_\al$, $\rL$}%
\index{$\rL_\al[u]$}%
\index{$\rL[u]$}%
\index{$\rL_\al$}%
\index{$\rL$}%
\eay
\right\}
}
\elin
%

\vyk{
The standard 
reference to $\ZF$ transfinite recursion 
does not work as it is based on the 
{\axf Collection/Replacement} schemata. 
Therefore we have to legitimize construction 
\eqref{LLdef} in $\zmc$ at a more basic level. 
}

\bte
[$\zmc$]
\lam{lae}
Suppose that\/ $u\sq\om$. 
Then
\ben
\renu
\itlb{lae1}%
each\/ $\rL_\al[u]$ is transitive and\/ 
$\al\sq\rL_{\al}[u]\,;$

\itlb{lae2}%
if\/ $\al<\ba$ then 
$\rL_\al[u]\in\rL_\ba[u]$ and\/ 
$\rL_\al[u]\sq\rL_\ba[u]\,;$

\itlb{lae3}%
if\/ $\la$ is limit then 
$\rL_\la[u]$ is closed under the 
rudimentary operations\/ 
\ref{rud} in Section~\ref{2prel}$\,;$

\vyk{
\itlb{lae2x}%
\sloppy
if\/ $n\ge2$ and\/ $X\sq\rL_\al[u]^n$ is definable over\/ 
$\rL_\al[u]$ by a formula with\/ $n$ free variables and 
parameters in\/ $\rL_\al[u]$ then\/ 
$X\in\rL_{\al+2n-1}[u]\,;$
}

\itlb{13c1}%
{\rm(I)}  
if\/ $\la\in\Ord$ is limit then the 
map\/ $\al\mto\rL_\al[u]$ $(\al<\la)$ is 
definable over\/ $\rL_\la[u]$ with $u$ as the 
only parameter, 
 {\rm(II)}  
the class-map\/ $\al\mto\rL_\al[u]$\/ 
$(\al\in\Ord)$  is 
definable over\/ $\rL[u]$ 
with $u$ as the only parameter$\,.$
\een
\ete

\bpf
See \cite{simp}, Theorem~VII.4.3 on 
\ref{lae1}, \ref{lae2}, \ref{lae3}. 
\vyk{
To prove \ref{lae2x} note that 
$\ans{x,y}\in\rL_{\al+1}[u]$ and 
$\ang{x,y}\in\rL_{\al+2}[u]$ for any $x,y\in\rL_\al[u]$, 
and hence $\ang{x_1,\dots,x_n}\in\rL_{\al+2n-2}[u]$ 
for $x_1,\dots,x_n\in\rL_\al[u]$ by induction, 
therefore\/ $\rL_\al[u]^n\sq\rL_{\al+2n-2}[u]$, and 
then easily 
$$
X=\ens{\ang{x_1,\dots,x_n}}
{x_1,\dots,x_n\in\rL_\al[u]\land  
\rL_{\al+2n-2}[u]\mo \vpi^{\rL_\al[u]}(x_1,\dots,x_n)}
$$
belongs to $\in\rL_{\al+2n-1}[u]$.   
By the way, 
using the $\al$-pair $\ang{x,y}_\al$, 
defined as in \cite[B.5, Lemma 6.1]{handb}, 
we get $\ang{x_1,\dots,x_n}_\al\in\rL_{\al}[u]$
for all $x_1,\dots,x_n\in\rL_\al[u]$. 
}%
Regarding \ref{13c1}, 
see Theorem~VII.4.8 in \cite{simp} or 
\cite[B.5, Lemma 4.1]{handb} in case $u = \pu$. 
\epf

What kind of set theory is provided in $\rL[u]$ by $\zmc$?

\ble
[$\zmc$]
\lam{nosep}
Let\/ $u\sq\om$.
All axioms of\/ $\zm$ sans\/ \Sep\  
hold in\/ $\rL[u]$ and in any set\/ $\rL_\la[u]$, 
where\/ $\la\in\Ord$ is limit.
\ele
\bpf[sketch]
This does not differ from the full-$\ZF$ case. 
Consider \eg\ the {\sf Union} axiom. 
Let $X\in\rL[u]$, so that $X\in\rL_\al[u]$, 
$\al\in\Ord$. 
As $\rL_\al[u]$ is transitive, 
the union $Y=\bigcup X\sq\rL_\al[u]$ 
is definable over $\rL_\al[u]$, 
hence $Y\in\rL_{\al+1}[u]=\Def{\,\rL_\al[u]}$. 
\epf

Axioms of $\zmc$ do not imply 
that the schemata of {\axf Replacement/Collection} 
necessarily 
hold in $\rL$, as the next example shows. 

\bpri
\lam{noco}
Arguing in the full $\ZF$ theory, 
let $\gM=\rL_{\vt}$, 
where $\vt=(\aleph_\om)^\rL.$ 
Let $\gN$ be the forcing extension of $\gM$ by 
ajoining a generic sequence of (generic) maps 
$f_n:\om\onto (\aleph_n)^\rL.$ 
Then $\gN$ is a model of $\zmc$. 
However $(\rL)^\gN=\gM$, and \Rep/\Col\ 
definitely fail in $\gM$. 
\epri

Unlike \Rep/\Col,
the {\sf Separation} schema always holds 
in $\rL$ under the $\zmc$ axioms in the 
background set universe by Theorem~\ref{3t}.

\parf{Definability and wellorderings} 
\las{abs}

Our goal here is to prove a few more delicate 
results related to the constructible hierarchy. 
The next lemma presents a key definability result. 

\vyk{
\bte
[$\zmc$]
\lam{defA}
Assume that\/ $u\sq\om$, and either\/ $K=\rL[u]$, 
or\/ $K=\rL_\la[u]$, where\/ $\la>\om$ is limit. 
Let\/ $Y\in K,$ 
$\al\in\Ord$, $\al<\la$ in case $K=\rL_\la[u]$, 
and let\/ 
$\vpi$ be a closed\/ $\in$-formula with parameters 
in $\rL_\al[u]$. 
Then the relation\/ $\rL_\al[u]\mo\vpi$ 
is absolute for $K.$   
\vyk{
\ben
\renu
\itlb{defA1}%
if\/ $X\in K$ is transitive and\/ $\vpi\in\form[X]$ 
then\/ $X\mo\vpi$ iff\/ $(X\mo\vpi)^K;$  
 
\itlb{defA2}%
if\/ $X,Y\in K$  
then\/ $Y=\Def X$ iff\/ $(Y=\Def X)^K.$ 
\een
}
\ete
\bpf
\vyk{
We'll have to use the {\em Quine--Rosser pair\/} 
\bce
$\pqr{A,B}=\ens{\vpi(a)}{a\in A}\cup
\ens{\vpi(b)\cup\ans0}{b\in B}$, where 
$\vpi(x)=(x\bez\om)\cup\ens{n+1}{n\in x\cap\om}$,
\ece
and the corresponding tuples 
$\pqr{A_1,\dots,A_n}$ for any $n\ge 3$, instead of the 
ordinary Kuratowski pairs and tuples. 
Their advantage is that if $A_1,\dots,A_n\in\rL_\al[u]$ 
then $\pqr{A_1,\dots,A_n}\in\rL_\al[u]$ as well, which 
holds for the standard tuple definition only in case 
$\al$ limit, which is not sufficient for this proof.
}%
In the context of the theorem, $\in$-formulas are 
viewed as tuples (finite sequences) of logical symbols 
and sets used as parameters, see \eg\ 
\cite[B.5, \S\,4]{handb}. 
Let $X=\rL_\al[u]$. 
Let $\form[X]\sq X\lom$ be the collection of all 
\imar{form[X]}%
\index{$\form[X]$}%
$\in$-formulas with parameters $p\in X.$ 
Thus if $\vpi\in\form[X]$ then $\vpi\in X^n,$ 
where $n=\lh\vpi$ is the length 
(the number of terms) of  
$\vpi$ as a tuple. 

Accordingly if $\vpi\in\form[X]$ then 
the set $\form_\vpi[X]$ of all subformulas $\psi$ 
of $\vpi$ (including $\vpi$ itself), 
in which all free variables are replaced 
by parameters $p\in X,$ satisfies 
$\form_\vpi[X]\sq X^n,$ where $n=\lh\vpi$.
\imar{form vpi[X]}%
\index{$\form_\vpi[X]$}%

Then the relation $X\mo\vpi$ is the formula 
$\Sat(X,\vpi)$ saying that $X$ is transitive, 
$\vpi\in\form[X]$ is closed, and there is a 
{\em valuation\/} 
$f=f_{X\vpi}:\form_\vpi[X]\to\ans{0,1}$, 
which satisfies the standard Tarski rules, and such 
that $f(\vpi)=\text{\sc true}:=1$. 
See \eg\ \cite[B.5, \S\,3]{handb}, 
\cite[p.~273]{simp}, or elsewhere. 

The notion of a valuation is obviously absolute 
for $K$. 
Therefore to prove that $X=\rL_\al[u]\mo\vpi$ is 
absolute it suffices to check that $f_{X\vpi}\in K$ 
for any (not necessarily closed) $\vpi\in\form[X]$.
This is proved by induction on the logical complexity 
of $\vpi$. 

For elementary formulas $\vpi$ of the form 
$v\in v'$ ($v,v'$ are free variables), we have  
$\form_\vpi[X]=\ens{\mlaq{x\in y}}{x,y\in X}$, 
and $f_{X\vpi}(\mlaq{x\in y})=1$ if in fact $x\in y$, 
and $=\,0$ otherwise, 
and $f_{X\vpi}\in K$ follows by a simple argument. 

As for the step, let say 
$\vpi=\vpi(\vec v)$ be $\psi(\vec v)\land\ta(\vec v)$, 
where $\vec v$ is a list $v_1,\dots,v_m$ of 
free variables.
By the inductive hypothesis, 
$f_{X\psi}$ and $f_{X\ta}$ belong to $K,$ hence 
to $\rL_\ga[u]$ for some $\ga\in\Ord$, 
$\ga<\la$ in case $K=\rL_\la[u]$. 
Then  
$\form_\vpi[X]=\form_\psi[X]\cup\form_\ta[X]\cup
\ens{\vpi(\vec x)}{\vec x\in X^m}.$ 
Accordingly $f_{X\vpi}$ is the union of 
$f_{X\psi}$ and $f_{X\ta}$ plus the   
assignment that $f_{X\vpi}(\vpi(\vec x))=1$ iff 
$f_{X\psi}(\psi(\vec x))=f_{X\ta}(\ta(\vec x))=1$ 
for all $\vec x\in X^m.$  
Thus $f_{X\vpi}$ is definable over $\rL_\ga[u]$ 
(with parameters $f_{X\psi},f_{X\ta}$), 
so $f_{X\vpi}\in \rL_{\ga+1}[u]\sq K$. 

The step of $\vpi(\vec v)$ being 
$\sus w\,\psi(\vec v,w)$, is 
considered the same way: $f_{X,\vpi(\vec v)}$ is equal 
to $f_{X,\psi(\vec v,w)}$,  
augmented by the assignment 
$f_{X\vpi(\vec v)}(\vpi(\vec x))=1$ iff 
$f_{X,\psi(\vec v,w)}(\psi(\vec x,y))=1$ for 
at least one $y\in X.$ 

This completes the proof of the theorem.
\epf 
}

\ble
[$\zmc$]
\lam{lae4}%
Let\/ $u\sq\om$, $\la$ be limit, and\/ 
$Y\in\rL_{\la}[u]$.  
Then\/ $Y$ is definable over\/ $\rL_\la[u]$ 
\/ {\rm(i)} 
by a formula with parameters\/ $\rL_\da[u]\yd\da<\la$, 
and\/ {\rm(ii)} 
by a formula with parameters\/ $\da<\la$ and\/ $u$. 
\ele

\bpf
(i) 
By definition, 
$Y=\ens{y\in\rL_\al[u]}{\rL_\al[u]\mo\vpi(y)}$, 
where $\al<\la$ and 
$\vpi$ may contain parameters in $\rL_\al[u]$. 
Arguing by induction on $\al$, 
let say $\vpi(y)$ be $\vpi(p,y)$, 
where $p\in\rL_\al[u]$ is a parameter. 
Then $p\in\rL_{\ga+1}[u]$ for some $\ga<\al$ 
by \eqref{LLdef} above. 
According to the inductive hypothesis, we have 
$p=\ens{z\in\rL_\ga[u]}{\rL_\la[u]\mo\psi(z)}$, 
where $\psi$ has only sets $\rL_\da[u]\yd\da<\la$, 
as parameters. 
Then 
$Y=\ens{y\in\rL_\al[u]}{\rL_\la[u]\mo\Phi(y)}$, 
where 
$$
\Phi(y) \, := \, 
\sus p\,
\big(y,p\in\rL_\al[u]\land 
p=\ens{z}{z\in\rL_\ga[u]\land \psi(z)}\land 
\vpi(p,y)^{\rL_\al[u]}
\big), 
$$ 
\vyk{
Consider the formulas 
\pagebreak[0]
\blin
\bul
{lae6}
{\left.
\bay{lcl}
\psi'(z) &:=& 
z\in\rL_\ga[u]\land \psi(z)  
, \\[0.7ex]
\vpi'(y) &:=& 
\sus p\,
\big(
p=\ens{z}{z\in\rL_\ga[u]\land \psi(z)}\land \vpi(p,y)
\big).%
\eay
\right\}
}
\elin
where
}%
and $\vpi(p,y)^{\rL_\al[u]}$ 
means the formal relativization 
to $\rL_\al[u]$, \ie, all quantifiers 
$\sus a,\,\kaz a$ are changed to resp.\  
$\sus a\in\rL_\al[u]\yd\kaz a\in\rL_\al[u]$. 
Then $\vpi'$ has only the sets $\rL_\ga[u]$, 
$\rL_\al[u]$, and some 
$\rL_\da[u]\yd\da<\la$, as parameters. 
%
This proves part (i). 
Now to infer part (ii) apply Theorem~\ref{lae}\ref{13c1}.
\epf

\ble
[$\zmc$]
\lam{newl}
Let\/ $u\sq\om$ and\/ 
$\la$ be limit. 
There is a map\/ 
$H:D=\om\ti\la\ti\la\lom$ onto $\rL_{\la}[u]$, 
definable over $\rL_\la[u]$ with\/ $u$ as 
the only parameter. 
\ele
\bpf
By Lemma~\ref{lae4}, 
each $Y\in\rL_\la[u]$ has the form 
$Y=\ens{y\in\rL_\al[u]}{\rL_\la[u]\mo\vpi(y)}$ 
for some $\al<\la$, 
where $\vpi$ contains parameters $\da<\la$ and $u$. 

Given a triple of $n,\al,p$ 
of $n\in\om$, $\al<\la$, and  
$p=\ang{\da_1,\dots,\da_k}\in\la^k$,  
let $\vpi_n$ be the $n$-th  \paf\  
$\in$-formula. 
If 
\ben
\fenu
\atc
\itlb{newl*}%
$\da_1\dots,\da_k<\la$ and 
$\vpi_n$ is $\vpi_n(v_1,\dots,v_k,v)$ 
with $k+1$ free variables,
\een
then define the set 
\bce
$H(n,\al,p)=\ens{y\in\rL_\al[u]}
{\rL_\la[u]\mo\vpi(\da_1,\dots,\da_k,y)}$, 
\ece
while if \ref{newl*} fails then just put $H(n,\al,p)=\pu$. 
Then $H$ 
is definable over $\rL_\la[u]$ with $u$ as a parameter 
by Theorem~\ref{lae}\ref{13c1} since it is defined 
in terms of the definable map $\al\mto\rL_\al[u]$. 
\epf

\vyk{
\bcor
[$\zmc$]
.\qed 
\ecor

\bcor
[$\zmc$]
\lam{13c2}
If\/ $\la\in\Ord$ is limit and\/ 
$Y\in\rL_{\la+1}$ 
then\/ $Y$ is definable in\/ $\rL_\la$  
by a formula  
with $u$ and ordinals\/ $\ga<\al$ 
as parameters. 
In particular, this is true for all\/ 
$Y\in\rL_\la$.\qed 
\ecor

\bpf
Use Theorem~\ref{lae}\ref{13c1} 
and Lemma~\ref{lae4}.
\epf
}

\ble
[$\zmc$]
\lam{13c3}
Let\/ $u\sq\om$.
There is a\/ 
\weo\/ $\lcu u$ of\/ $\rL[u]$ 
definable over\/ $\rL[u]$ 
with\/ $u$ as the only parameter.  
\imar{lcu}%
\index{$\lcu u$, $\lclu\la u$, well-orderings}%
\index{well-orderings $\lcu u$, $\lclu\la u$}%
If\/ $\la\in\Ord$ is limit then there is a 
\weo\/ $\lclu\la u$ of\/ $\rL_{\la}[u]$ 
\imar{lclu la u}%
definable over\/ $\rL_\la[u]$ 
with\/ $u$ as the only parameter.  
\ele

\bpf
In the $\la$-case, let the map 
$H:D\onto \rL_{\la}[u]$ 
be given by Lemma~\ref{newl}. 
The set $D=\om\ti\la\ti\la\lom\sq \rL_\la[u]$ 
is \paf\ definable over $\rL_\la[u]$. 
Thus to define $\lclu\la u$ 
it suffices to show that  $D$ admits a \weo\ 
$<_D$ \paf\ definable over $\rL_\la[u]$. 
For that purpose, if 
\bce
$d=\ang{n,\al,u=\ang{\ga_1,\dots,\ga_m}}\in D$, \ 
 \  
$d'=\ang{n',\al',
u'=\ang{\ga'_1,\dots,\ga'_{m'}}}\in D$, 
\ece 
then 
let $\mu(d)=\max\ans{\al,\ga_1,\dots,\ga_m}$ 
and define $d<_D d'$ if and only if: 
\ben
\fenu
\atc\atc
\itlb{13c4}%
either $\mu(d)<\mu(d')$,\\  
or $\mu(d)=\mu(d')$ and $m<m'$,\\  
or $\mu(d)=\mu(d')$, $m=m'$, and $u<u'$ 
lexicographically in $\la^m$,\\  
or $\mu(d)=\mu(d')$, $m=m'$, $u=u'$, and  
$n<n'$.  
\een

The \weo\ $\lcu u$ of $\rL[u]$ is then defined so that 
$x\lcu u y$ iff either (1) $\la_x<\la_y$, where $\la_x$ 
is the least limit ordinal with $x\in\rL_{\la_x}$, 
or (2) $\la_x=\la_y$ and $x\lclu\la u y$. 
\epf

\vyk{
\bcor
[$\zmc$]
\lam{13c5}
Let\/ $u\sq\om$.
If\/ $\la\in\Ord$ is limit\/ \kn{not} of the 
form\/ $\al+\om$, then every set\/ 
$X\in \rL_{\la}[u]$ 
admits a \weo\ $\prec$ in\/ $\rL_\la[u]$.  
\ecor

\bpf
Let $X\in\rL_\al[u]$, $\al<\la$. 
The restriction ${\prec}={{\lclu{\la+\om}u}\res X}$ 
is definable over $\rL_{\al+\om}[u]$ with $X$ as 
a parameter, by Lemma~\ref{13c3}, 
hence ${\prec}$ belongs to 
$\rL_{\al+\om+1}[u]\sq\rL_\la[u]$.
\epf
}

\vyk{\gol
The next result is interesting in the 
context of the discussion in the end of 
Section~\ref{2r}.

\bcor
[$\zmc$, essentially Lemma 1 in \cite{BPrah}]
\lam{13c5}
Let\/ 
$\ens{k\in\om}{\vpi^\rL(k)}\nin\rL$, where\/  
$\vpi(k)$ is an\/ $\in$-formula with parameters 
in\/ $\rL$. 
Then there is a \paf\ formula\/ $\psi(k)$ such 
that still\/   
$\ens{k\in\om}{\psi^\rL(k)}\nin\rL$.
\ecor

\bpf
Let $\vpi(k)$ be say $\vpi(p_0,k)$, where 
$p_0\in\rL$ is the only parameter. 
Let 
\bce
$\Phi(p):=\mlap{
\ens{k\in\om}{\vpi(p,k)}\text{ is not a set}}$, 
\ece
so that, under our assumptions, we have 
$\sus p\in\rL\,\Phi^\rL(p)$. 
Let $p_1$ be the $\lc$-least $p\in\rL$ 
satisfying $\Phi^\rL(p_1)$. 
Then $p_1$ is \paf\ definable over $\rL$ since 
so is $\lc$ by Lemma~\ref{13c3}.
Thus $\psi(k):=\Phi(p_1,k)$ is as required.
\epf
}

\parf{The master plan}
\las{mas}

The purpose of this Section is to formulate a 
convenient necessary condition for getting $\zfcm$ 
in some constructible domains 
(Theorem~\ref{z37mz}). 
%
%
%
To simplify formalities, define the following 
formula:

\bdf
[$\zmc$]
\lam{sig}
Let $\ssu(u,\Om,K)$ say: 
$u\sq\om$, and\vim
\index{$\ssu(u,\Om,K)$}%
\imar{\ \ \ , ssu}%
\ben
\item[$-$]
either 
(A) $\Om=\Ord$, $K=\rL[u]$, and $\omiu u$ 
does not exist, in other words, every ordinal 
is countable in $\rL[u]$,\vim

\item[$-$]
or 
(B) $\Om=\omiu u$ exists, 
and $K=\rL_\Om[u]=\rL_{\omiu u}[u]$.\vim
\een

\noi
Thus 
$K=\bigcup_{\al\in\Om}\rL_\al[u]$ in both cases (A), (B).  
\vyk{
\itlb{sig3}%
\ocoun\ holds in $K$, 
\itlb{sig4}%
in case (B), there is no maps\/ 
$G:\om\text{ onto }\Om\yd G\in\rL[u]$. 
Thus $K=\rL[u]$ in case (A) and 
$K=\rL_\Om[u]=\rL_{\omiu u}[u]$ in case (B). 
%
\index{$<_K$}%
}
\edf

\ble
[$\zmc+\ssu(u,\Om,K)$]
\lam{2coun}
If $\al\in\Om$ then $\rL_\al[u]$ 
is ctble in $\rL[u].$ 
\ele
\bpf
Let $\al\in\Om$ be limit. 
By Definition~\ref{sig}, there is a map $f\in \rL[u]$,
$f:\om\text{ onto }\al$. 
Lemma~\ref{lae4} provides a set 
$D=\om\ti\al\ti\al\lom\in \rL[u]$ and a map 
$H\in\rL[u]\yd H:D\text{ onto }\rL_{\al}[u]$. 
%
We get a map 
$h\in\rL[u]\yd h:\om\text{ onto }\rL_\al[u]$, 
by combining $f$ and $H$ in $\rL[u]$. 
\epf

\ble
[$\zmc+\ssu(u,\Om,K)$]
\lam{z33mz}
Assume that\/ $X\in K$, and\/ 
$F:X\to K$ is a class-function definable 
over\/ $\rL[u]$. 
Then\/ $\ran F=\ens{F(x)}{x\in X}\sq\rL_\ga[u]$ 
for some\/ $\ga\in\Om$, 
hence\/ 
$F\yd\ran F$ are sets. 
\ele

\bpf
By Lemma~\ref{2coun}, 
we \noo\ suppose that $X=\om$. 
For any $k<\om$, let $\da_k$ be the least $\da\in\Om$ 
satisfying $F(k)\in \rL_\da[u]$. 
Assume towards the contrary that $\ens{\da_k}{k<\om}$ 
is {\kn {\rm un}}bounded in $\Om$. 
Then $\Om=\bigcup_{k<\om} {\da_k}$. 

In case (A), 
for any $k$, let $h_k$ be the $\lcu u$-least
of all functions $h\in \rL[u]$, 
$h:\om\text{ onto }\da_k$. 
(Such $h$ exist by Definition~\ref{sig}.)  
If $n = 2^k(2j+1)-1$  
then put $G(n)=h_k(j)$. 
Then $G$ is a definable class-function 
from $\om\text{ onto }\Om=\Ord$ by construction. 
Thus $\Om$ and $G$ are sets 
by Lemma~\ref{mz12} since $\Om$ is transitive. 
This is a contradiction  
since $\Ord$ is not a set in $\zmc$. 
 
In case (B), $\Om=\omiu u$. 
Define $h_k$ and then $G$ using the \weo\/ 
$\lclu\Om u$ of\/ $\rL_{\Om}[u]$ instead of $\lcu u$. 
Then $G$ is a class-function 
from $\om\text{ onto }\Om=\omiu u$, definable over 
$\rL_\Om$ since such is $\lclu\Om u$. 
Thus $G\in\rL_{\Om+1}[u]\sq\rL[u]$, and hence 
$\Om$ is countable in $\rL[u]$. 
This is a contradiction.  
\epf 

\vyk{
Thus $\Om$ is a limit ordinal and  
$K=\rL_\Om$ by $\ssu(u,\Om,K)$. 
Therefore $G\in\rL_{\Om+1}[u]$ by construction.
This contradicts \ref{sig4} of Definition \ref{sig}.
}

\bcor
[$\zmc+\ssu(u,\Om,K)$]
\lam{z35mz}
\sloppy
Assume that\/ $\al\in\Om$, $m<\om$, and\/ 
$G_1,\dots,G_m:K\to K$ be class-functions 
definable over\/ $\rL[u]$. 
There is a limit ordinal\/ 
$\ba\in\Om\yi\ba>\al$, satisfying\/ 
$\imb{G_k}{\,\rL_\ba[u]}\sq\rL_\ba[u]$ for all\/ 
$k=1,\dots,m$.
\ecor

\bpf
Put $G(x)=\ang{G_1(x),\dots,G_m(x)}$.
Use Lemma~\ref{z33mz} to get a class-sequence 
$\al=\al_0<\al_1<\al_2<\dots$ of ordinals in 
$\Om$ satisfying 
$\imb G{\rL_{\al_n}[u]}\sq\rL_{\al_{n+1}}[u]$, 
$\kaz n$. 
Apply Lemma~\ref{z33mz} again to show that 
$\ba=\sup_n\al_n\in\Om$. 
\epf

Assume $\ssu(u,\Om,K)$.
Say that 
$\ba\in\Om$ {\em reflects} a 
\index{reflects}%
formula $\vpi(x_1,\dots,x_n),$  
if 
the equivalence
$\vpi^{K}(x_1,\dots,x_n) \eqv 
\vpi^{\rL_\ba[u]}(x_1,\dots,x_n)$
holds for all $x_j\in\rL_\ba$.
The following {\em reflection lemma} 
is a standard consequence of 
Corollary~\ref{z35mz}. 

\ble
[$\zmc+\ssu(\Om,K)$]
\lam{z36mz}
If\/ $\al\in\Om$ and\/ $\vpi$ is 
a \paf\ formula
then there exists a limit ordinal\/ 
$\ba\in\Om\yi\ba>\al$ 
which reflects\/ $\vpi$ and every subformula 
of\/ $\vpi$. 
\ele
\bpf[sketch]
We \noo\ assume that $\vpi$ does not contain $\kaz$ (otherwise replace $\kaz$ with $\neg\,\sus\neg$). 
Let's enumerate $\psi_1,\dots,\psi_n$ 
all the sub-formulas of $\vpi$ 
(including possibly $\vpi$ itself) 
beginning with $\sus$. 
If $j\le n$ then we define a class-function 
$G_j$ as follows. 

Assume that $j\le n$ and   
$\psi_j$ is $\sus y\,\chi_j(y,x_1,\dots,x_m)$.  
If $p=\ang{x_1,\dots,x_m}\in K$ and 
there is $y\in K$ satisfying 
$\chi_j^{K}(y,x_1,\dots,x_m)$,  
then let $G_j(p)$ 
be the $\lcu u$-least of these $y$. 
Otherwise let $G_j(p)=\pu$. 
Each class-function $G_j$ 
is definable over $\rL[u]$ since such is 
the \weo\ $\lcu u$.  

\sloppy
By Corollary~\ref{z35mz}, there is an ordinal 
$\ba\in\Om$, $\ba>\al$,
satisfying\/ 
$\imb{G_j}{\,\rL_\ba[u]}\sq\rL_\ba[u]$ 
for all\/ $j=1,\dots,n$.
Now it easily goes by induction on the number 
of logical symbols  
that $\ba$ reflects every subformula of $\vpi$, 
in particular it reflects $\vpi$ itself, 
as required.
\epf

\bte
[$\zmc+\ssu(u,\Om,K)$]
\lam{z37mz}
\Sep\ and \Col\ 
hold in\/ $K$. 

Therefore\/ $\zfcm$ as a whole  
holds in\/ $K$  by Lemma~\ref{nosep}.
\ete

\bpf
{\axf Separation}. 
Assume that 
$\vpi(x,y)$ is a parameter-free formula, 
$\al\in\Om$, $p\in X=\rL_\al[u]$. 
We have to prove that  
$Y=\ens{x\in X}{\vpi^{K}(x,p)}\in K$. 
Let, by Lemma~\ref{z36mz}, 
a limit ordinal $\ba\in\Om\yi\ba>\al$ 
reflect\/ $\vpi(x,y),$  so that 
\bce
$Y=\ens{x\in X}{\vpi^{\rL_\ba[u]}(x,p)} 
=\ens{x\in X}{{\rL_\ba[u]}\mo\vpi(x,p)} 
\in\rL_{\ba+1}[u]\sq K$. 
\ece
 
{\axf Collection}. 
Assume that   
$\vpi(x,y,z)$ is a parameter-free formula, 
$\al\in\Om$, $p\in X=\rL_\al[u]$, 
and we have 
$\kaz x\in X\,\sus y\in K\,
\vpi^{K}(x,y,p)$. 
By Lemma~\ref{z36mz}, there exists 
a limit ordinal $\ba\in\Om\yi\ba>\al$ 
which reflects\/ $\sus y\,\vpi(x,y,z)$, 
with all its subformulas, 
including $\vpi(x,y,z)$, so that 
$$
\hspace*{-1ex}
\kaz x\in X\,\sus y\in\rL_\ba[u]\,
\vpi^{\rL_\ba[u]}(x,y,p), 
\text{ and }\:
\kaz x\in X\,\sus y\in\rL_\ba[u]\,
\vpi^{K}(x,y,p).\hspace*{-4ex} 
\eqno\qed
$$
\ePf

\vyk{
\bte
[$\zmc$]
\label{2red}
\ben
\cenu
\itlb{2red2}%
$\rL\kpo$ 
satisfies\/ $\pad$. 

\itlb{2red1}%
$\Lz=\bigcup_{\al\in\Ord\text{ countable in }\rL}\rL_\al$
satisfies\/ $\zfcm.$ 
We may note here that
\ben
\raenv
\itlb{2red2A}%
if\/ $\omil=\Om$ exists then\/ $\Lz=\rL_\Om\,;$

\itlb{2red2B}%
if\/ $\omil$ does not exist then\/ $\Lz=\rL\,.$
\een
\een
\ete
}

\parf{Proof of Theorems~\ref{2red} and \ref{1t}}
\las{fura1}


Theorem~\ref{1t} is an elementary consequence 
of Theorem~\ref{2red}, so we concentrate on the 
latter.\vom

{\ubf Case (b) of Theorem~\ref{2red}.} 
Arguing in $\zmc$, we have case (B) of 
Definition~\ref{sig}
with $u=\pu$, $\Om=\omil$, $K=\Lz=\rL_{\omil}$. 
Then $\ssu(\pu,\omil,\Lz)$ holds, and 
hence $\Lz$ satisfies $\zfcm$ by 
Theorem~\ref{z37mz}.\vom

{\ubf Case (a) of Theorem~\ref{2red}.}
Similar, but via case (A) of Definition~\ref{sig}.


%
\vyk{
\bpf[Claim \ref{2red1a} of Theorem~\ref{2red}]
Still arguing in $\zmc$, 
we assume that \ocoun\ \rit{fails} in $\rL$, so that 
there are ordinals uncountable in $\rL$. 
In this case, the least one of them is 
denoted by $\omil$ in set-theoretic studies. 

It is pretty natural now to prove that 
$\ssu(\pu,\omil,K)$ holds for 
$K=\rL_{\omil}=\bigcup_{\al<\omil}\rL_\al$, and hence 
$K$ satisfies $\zfcm$ by Theorem~\ref{z37mz}, finally 
getting Claim \ref{2red1b} of Theorem~\ref{2red} with 
$\la=\omil$. 
However, going this way, we have to prove that 
\ocoun\ holds in $K=\rL_{\omil}$. 
This is doable, and in fact a similar claim is 
established in \cite{simp} in the course of the proof 
of Theorem VII.4.34. 
However the proof there involves rather special 
reasoning based on the $\is{}1$-theory of constructible 
hierarchy, the G\"odel collapse argument, and some other 
intermediate $\zfc$-style claims, 
which we don't plan to be used in our proof. 
Fortunately, there is a much simpler way 
to the same goal (\ref{2red1a} of Theorem~\ref{2red}), 
which gives a somewhat less sharp final result. 

If $u\sq\om$ then let $\xom[u]$ denote the least 
limit ordinal $\Theta$ such that 
\blin
\bul
{xom}
{\neg\:\sus f\in\rL_{\xom[u]+1}[u]\,
(f:\om\text{ onto }\xom[u]),
}%
\elin
if such limit ordinals $\Theta$ exist, and 
$\xom[u]=\Ord$ otherwise. 
Let $\rlau u=\bigcup_{\al\in\xom[u]}\rL_\al[u]$. 
%
\index{$\rlau u$}%
\index{class $\rlau u$}%
\index{set!constructible, $\rlau u$}%
\index{ordinal!$\xom[u]$}%
Let $\meb u$ be the \weo\ 
$\lclu{\xom[u]}u$ in case $\xom[u]\in\Ord$, 
and be the \weo\ $\lcu u$ in case $\xom[u]=\Ord$ 
(see Lemma~\ref{13c3}). 
Thus $\meb u$ is a definable 
\weo\ of $\rlau u$ in each case.
\index{wellorda@\weo\ $\meb u$}%
\index{$\meb u$, \weo}%

Rather similar ideas 
(for instance, the notion of a gap ordinal) 
are known in the studies of constructibility 
since 1970s, see  
\cite{mz81} (our source of the construction based on 
\eqref{xom}), 
and \eg\ \cite{marekA,marekG}. 


\ble
[$\zmc$, pretty obvious]
\lam{furs1}
Let $u\sq\om$. 
Then $\ssu(u,\xom[u],\rlau u)$ holds. 
Therefore\/ $\rlau u$ satisfies\/ $\zfcm$ 
by Theorem~\ref{z37mz}.\qed
\ele

Now to prove  
Claim \ref{2red1a} of Theorem~\ref{2red} 
apply Lemma~\ref{furs1} for $u=\pu$.
\epf
}

\vyk{

{\gol
\bpf[Theorem~\ref{2red}\ref{2red1}]
Arguing in $\zmc$, 
assume that $\Om=\omil$ exists. 
Let $K=\rL_{\Om}=\bigcup_{\al<\Om}\rL_\al$. 
By Theorem~\ref{z37mz}, it suffices to prove 
$\ssu(\Om,K)$.
\epf
}

\bpf[Theorem~\ref{1t}]
Make use of Corollary~\ref{red1} and 
Theorem~\ref{2red}\ref{2red1}.
\epf

{\gol
{\ubf Simpson \cite[VII.4, VII.5]{simp}.}

$\HCL[u]$ is the collection of all sets $x$ such that 
there exists a map $f\in\rL[u]$ such that 
$\dom f=x$, $\ran f$ is a transitive set, and $x\in\ran f$. 

VII.4.25. Let $u\sq\om$. 
Then $\HCL[u]$ is an elementary submodel of $\rL[u]$ 
\poo\ all $\is{}1$ $\in$-formulas. 
By the G\"odel collapse technique. 

Corollary: $\HCL[u]\mo$ every set $x$ is constructible 
in $u$, 
that is, $\sus \al\in\Ord\,(x\in\rL_\al[u])$. 

Corollary. If $\Ord\sq\HCL[u]$ then 
$\HCL[u]=\rL[u]$. 

VII.4.34. 
If $\Ord\sneq\HCL[u]$ then 
$\HCL[u]=\rL_\vt[u]$, where $\vt={\omi^{\rL[u]}}$. 
}
}%
%
%


\parf{Proof of Theorem~\ref{3t}\ref{3tA}}
\las{fura2}

Arguing in $\zmc$, we are going to prove that 
$\rL\kpo$ satisfies\/ $\pad$. 
The argument breaks down into two cases.\vom

{\sl Case 1\/}: 
there is $u\sq\om$ such that $\omiu u$ does not exist.
Then $\zfcm$ holds in $\rL[u]$ 
by Theorem~\ref{z37mz}.  
Therefore $\zfcm$ holds in $\rL$ as well, and this 
implies $\pad$ in $\rL\kpo$, as required.\vom

{\sl Case 2\/}: 
$\omiu u\in \Ord$ exists for all $u\sq\om$. 
In particular, $\Om=\omil\in\Ord$ exists, and 
$\rL_\Om$ is a model of $\zfcm$ by 
Theorem~\ref{2red}  
already proved. 
Therefore it suffices to prove that $\rL\kpo\sq\rL_\Om$. 

This is a well-known result in $\zfc$ and $\zfcm,$ 
a part of G\"odel's proof of {\ubf CH} in $\rL$. 
G\"odel's reasoning is doable in $\zmc$,  
and a close claim is established in \cite{simp} 
in the course of the proof of Theorem VII.4.34. 
However the proof there involves quite special 
arguments, \eg\ the $\is{}1$-theory of constructible 
hierarchy, 
which we don't plan to use in our proof. 
Yet there is a much simpler way 
to the same goal, by reduction to the $\zfcm$ 
environment. 

Thus, to prove $\rL\kpo\sq\rL_\Om$, let $x\in\rL\kpo$. 
Then $x\in\rL_\la$ for some $\la\in\Ord$. %
We assert that
\ben
\fenu
\itlb{bcov}%
there is an ordinal\/ $\vt>\la$ such that\/ 
$\rL_\vt$ is a model of $\zfcm.$  
\een
Indeed, 
by \coun\ of $\zmc$, there is a bijection 
$h:\om\text{ onto }\la$. 
Let $u=\ens{2^j\cdot 3^k}{h(j)<h(k)}$ code $h$. 
Note that $\vt=\omiu u\in\Ord$ by the Case 2 
assumption, and $\rL_{\vt}[u]$ is a model 
of $\zfcm$ by Theorem~\ref{z37mz}, hence, 
$\rL_\vt\mo\zfcm$ as well. 
Thus it suffices to show that $\la\le\vt$.

Suppose to the contrary that $\vt<\la$. 
Then $\rL_{\vt}[u]\mo\zfcm$, as above. 
In addition, 
$\rL_{\vt}[u]$ is a model of ordinal height $\vt$,  
and $u\in\rL_\om[u]\sq\rL_{\vt}[u]$, by construction.
But $u$ effectively codes the ordinal
$\la>\vt$, which is a contradiction. 
This completes the proof of \ref{bcov}.  

Choose $\vt$ by \ref{bcov}; thus $x\in\rL_\vt$. 
We don't claim that $\Om=\omi^{\rL_\vt}$, 
but still $\Om$ obviously remains a regular 
uncountable cardinal in $\rL_\vt$. 
This implies, even in $\zfcm,$ that 
$\rL_\vt\kpo\sq\rL_\Om$ by a standard G\"odel 
collapse argument. 
We conclude that $x\in\rL_\Om$, 
as required.

\parf{Proof of Theorem~\ref{3t}\ref{3tB}, sketch}
\las{fura3}

We'll make use of a deep result 
in \cite{BPrah} related to countable 
\rit{index ordinals}, \ie, those $\al$ satisfying 
$(\rL_{\al+1}\setminus\rL_\al)\cap\pws\om\ne\pu$. 
{\ubf We argue in $\zmc$.}

First of all, it is asserted in \cite{BPrah} that 
there eexists a parameter-free closed 
$\in$-formula $\sg$ 
such that, for any transitive set $M$, 
$\sg^M$ (the formal relativization) holds 
iff $M=\rL_\la$ for some limit ordinal $\la$. 
Basically, $\sg$ says, in some proper way, 
that all sets are constructible and there is 
no maximal ordinal. 
The required property is based on the 
absoluteness of the G\"odel 
construction for all transitive sets satisfying some 
simple conditions, \cite{jechmill}.

Now, suppose to the contrary that \Sep\ fails 
in $\rL$, that is, there exist: a transitive 
set $B\in\rL$ (say $B=\rL_\al$ for some $\al$) and 
a formula $\vpi(p,x)$ with a parameter $p\in\rL$, 
such that $Y=\ens{b\in B}{\vpi^\rL(p,b)}\nin \rL$. 
($Y$ is a set in the 
background $\zmc$ universe by  \Sep\ of $\zmc$.) 
Taking the $\lc$-least $B$ and $p$ with these 
properties, we reduce the general case to the 
following: 
\ben
\fenu
\atc
\itlb{ls*}\msur%
$B=\ens{b\in\rL}{\vt^\rL(b)}$ 
is parameter-free definable in $\rL$, and 
$\Phi(x)$ is a parameter-free formula, still 
satisfying $Y=\ens{b\in B}{\vpi^\rL(b)}\nin \rL$.
\een

Let $f_1,\dots,f_m$ be the list of Skolem 
functions for all existential subformulas of the 
formulas 
\ben
\fenu
\atc
\atc
\itlb{ls**}\msur%
$\sg\yi\vpi(x)$, 
and `$B=\ens{b\in\rL}{\psi(b)}$', 
and their negations, 
\een
defined in $\rL$ on the basis of the 
parameter-free definable \weo\ $\lc$. 

Consider the closure $M$ of $B\cup\ans B$ under 
$f_1,\dots,f_m$. 
By a standard combinatorial argument, there is 
a 
definable class-map $\Phi$ defined on the set 
$U=B\lom\ti\om\lom$, such that 
$M=\Phi\ima U$. 
Let $\ta:M$ onto a transitive class $N$ be a 
collapse map, \ie, 
$\ta(x)=\ens{\ta(y)}{y\in x\cap M}$ for all 
$x\in M.$ 
(To define $N,\ta$ apply 
Corollary~\ref{c12} for sets $M_\al=M\cap\rL_\al$, 
$\al\in\Ord$, and let $\ta$ be the union 
of all partial collapse maps $\ta_\al:M_\al$ 
onto a transitive set $N_\al$.) 

Using Lemma~\ref{mz12} for the superposition of 
$\Phi$ and $\ta$, we conclude that $N$ is a set. 
Moreover, as $B$ is transitive, we have 
$B=\ta(B)\in N$. 

On the other hand, the class or set $M$ is an 
elementary submodel of $\rL$ \poo\ formulas 
\ref{ls**}  
by construction. 
In particular, $M\mo\sg$, hence $N\mo\sg$ as well, 
and we conclude by the choice of $\sg$ 
that $N=\rL_\la$ for some limit $\la$. 

By the same argument (and because $B=\ta(B)$) 
we conclude that 
$Y=\ens{b\in B}{\vpi^{\rL_\la}(b)}\in\rL_{\la+1}
\sq \rL$, 
which contradicts \ref{ls*}.

\parf{A corollary in the domain of reals}
\las{2r}

Theorem~\ref{2red} just proved implies the following 
corollary. 

\bcor
[$\padm$]
\lam{2redC}
$\rL\kpo$ satisfies\/ $\pad$. 
Saying it differently, $\rL\kpo$ is an 
interpretation of\/ $\pad$ in\/ $\padm$. 
\ecor

Here $\rL\kpo$ essentially 
means $\ens{x\sq\om}{\kons(x)}$, where $\kons(x)$ 
is a certain $\is12$ formula of $\lpad$ that 
expresses the constructibility of $x\sq\om$ by 
referring to the existence of a real that encodes 
(similar to \eg\ encoding by trees in $\wft$) 
a set-theoretic structure that
indicates the constructibility of $x$. 
Such a formula was explicitly defined by 
Addison~\cite{add2,add1}, but implicitly can be 
found in studies by G\"odel~\cite{godel40} and 
Novikov~\cite{nov1951}.

\bpf[sketch]
The $\padm$ structure $\bI$ satisfies $\zmc$ by 
Theorem~\ref{it1}. 
Therefore we have 
$\jnt{\text{$\rL\kpo$ satisfies\/ $\pad$}}$ 
by Theorem~\ref{2red}. 
Yet the $\bI$-reals are isomorphic to the 
true reals in the background $\padm$ universe.
We conclude that, in $\padm$, 
$\rL\kpo$ satisfies\/ $\pad$.
\epf

Corollary~\ref{2redC} can be compared with 
its better-known version:

\bpro
[Theorem~1.5 in \cite{marekA}]
\lam{2r*}
If $X\sq\pws\om$ 
is a $\ba$-model of\/ $\padm$ then $X\cap\rL$ 
is a $\ba$-model of\/ $\pad$ plus constructibility. 
\qed 
\epro

The proof in \cite{marekA} involves Lemma~1.4 
there, that cites Theorem~1 in \cite{BPrah}, 
presented in Proposition \ref{BP}\ref{BP1} below.
\vyk{
saying that if 
$\rL_{\al+1}\bez\rL_\al$ contains a set $x\sq\om$ 
then there is a relation $E\in\rL_{\al+1}$ on $\om$ 
isomorphic to the structure $\stk{\rL_\al}{\in}$. 
The proof in \cite{BPrah} in turn 
refers to some results in Boolos' PhD Thesis 
(unpublished). 
}%
Another path to Proposition~\ref{2r*}, 
quite complicated in its own way,  
is given in 
\cite{EF,ender}. 

It is definitely tempting to accomodate these 
proofs of Proposition~\ref{2r*} to the case  
$X=\pws\om$ towards Corollary~\ref{2redC} 
under the $\zmc$ axioms. 
Yet we are not going to pursue this plan 
here as it will definitely involve   
more complex  
arguments than the above proof of 
Theorems~\ref{2red} and \ref{3t}.


\vyk{
{\ubf The idea} of the proof of Theorem~\ref{red2} 
can be roughly presented as follows. 

Arguing in $\zmc$, 
we define the transitive class 
$\rL=\bigcup_{\al\in\Ord}\rL_\al$ 
of all G\"odel-constructible sets, 
with the hope that 
it satisfies $\zml$. 
We may note that the $\zmc$ ordinals are 
no less rich than the $\zfm$ ordinals, in 
particular, any countable \weo\ is isomorphic 
to an ordinal by \moco\ of Definition \ref{zmad}. 
Thus we may hope that the constructible sets 
behave in $\zmc$ not much worse than in $\zfm,$ 
leading to $\zml$ in $\rL$.  
However, what can lead to a failure here is 
\coun, 
which is not necessarily preserved in 
$\rL$ under the $\zmc$ axioms.
Therefore our plan will be 
to prove that either 
the whole $\rL$ or one of 
$\rL_\al$ ($\al$ limit) satisfies $\zml$. 
}

\parf{Some other models}
\las{2om}

Here we briefly describe three other models of 
$\zfcm$ in $\zmc$ which work similar to $\Lz$ 
of Theorem~\ref{2red}.\vom

{\ubf Model 1.} 
Consider the least ordinal $\La$ such that 
the set 
$\rL_\La$ is not countable in $\rL_{\La+1}$ 
--- provided such ordinals 
exist, and otherwise $\La=$ all ordinals.  
Put $\Ld=\bigcup_{\al\in\La}\rL_\al$. 
It is demonstrated in \cite{mz81} that 
$\Ld$ is an interpretation of $\zfcm$ 
in  $\zmc$. \vom

{\ubf Model 2: a version of Model 1.} 
Consider the least ordinal $\Xi$ such that 
the difference $\rL_{\Xi+1}\bez \rL_\Xi$ 
contains no sets $x\sq\om$ 
--- 
the first \rit{index ordinal} as defined in 
\cite{BPrah} 
--- provided such ordinals 
exist, and otherwise $\Xi=$ all ordinals.  
Arguments close to those in \cite{mz81} 
show that $\Ldd=\bigcup_{\al\in\Xi}\rL_\Xi$ 
$\Ld$ is an interpretation of $\zfcm$ 
in  $\zmc$.\vom

{\ubf Model\/ 3.} 
Simpson defines in \cite[VII.4.22]{simp} the 
set or class $\HCL$ 
of all sets $x$ which belong to 
transitive sets $X\in\rL$, countable in $\rL$, 
and proves that $\HCL$ is an interpretation 
of $\zfcm$ in  $\zmc$ yet again. 
But it looks like $\HCL$ is just equal to 
$\Lz$ of Theorem~\ref{2red}.


\parf{Ramified analytical hierarchy --- 
a shortcut\,?}
\las{rah}

Cutting Theorem~\ref{1t} to the equiconsistency 
of $\pad$ and $\padm$ (second order arithmetic 
with, resp., without the countable Choice $\ACw$), 
one may want to manufacture a true 
second-order arithmetical proof, not 
involving set theories like 
$\zm,\,\zfcm,\,\zfm,\,\zmc$. 
The above proof (Section~\ref{fura1}) 
definitely does not 
belong to this type, 
since it involves $\zmc$ in quite significant 
way.
In this section, 
we survey a possible approach to this 
problem.

Using earlier ideas of Kleene~\cite{klee} and 
Cohen~\cite{cohMM}, a transfinite 
sequence of countable sets $\rA_\al\sq\pws\om$ 
is defined in \eg\ \cite[\S\,3]{BPrah} 
by induction so that
\pagebreak[0] 
\blin
\bul
{Adef}
{\left.
\bay{lcl}
\rA_0 &=& \pfin\om= 
{\text{all finite sets }x\sq\om}\\[0.7ex]
%
%
\rA_{\al+1} &=& 
\Def{\,\rA_\al}
\text{ \ for all \ }\al\\[0.7ex]
\rA_\la &=& 
\bigcup_{\al<\la}\rA_\al \ \ 
\text{for all limit }\la\\[0.7ex]%
\rA &=& 
\bigcup_{\al\in\Ord}\rA_\al \,= \,
\text{all {\em ramified analytic\/} sets}
\index{$\rA_\al$, $\rA$, ramified analytic}%
\index{set!ramified analytic, $\rA_\al$, $\rA$}%
\index{ramified analytic!set}%
\index{ramified analytic!hierarchy}%
\eay
\right\},
}
\elin
where 
$\Def \rA_\al= \ens{x\sq\om}
{x\text{ is definable over }\rA_\al 
\text{ with parameters}}$ 
in the 2nd line. 
Thus a set $x\sq\om$ belongs to  
$\Def {\rA_\al}$ iff 
$x=\ens{n}{\rA_\al\mo\vpi(n)}$ for some 
formula $\vpi$ of $\lpad$ with parameters in 
$\rA_\al$, and $X\mo\dots$ means 
the formal truth in the $\lpad$-structure 
$\stk\om X$. 
The following is routine.

\ble
\lam{bhp}
If\/ $x\in \rA_\al$ and\/ $y\sq\om$ is 
arithmetical in\/ $x$ then\/ $y\in \rA_\al$.\qed
\ele 

In spite of obvious similarities with the 
G\"odel constructible hierarchy \eqref{LLdef}, 
the ramified analytic hierarchy is 
collapsing below $\omi$:

\ble
[Cohen] 
\lam{eb0}
There is an ordinal\/ $\bbo<\omil$ such that\/ 
$\rA_{\bbo}= \rA_{\bbo+1}=\rA_{\ga}$ for all\/ 
$\ga>\bbo$. 
\imar{bbo}%
\index{$\bbo$, ordinal}%
\index{ordinal $\bbo$}%
Then obviously\/ $\rA=\rA_{\bbo}$ and\/ 
$\rA\mo\padm$.
\ele

\bpf
By the cardinality argument, there is an ordinal 
$\ba$ with $\rA_{\ba}= \rA_{\ba+1}$. 
Then $\rA_\ba\mo\text{\Sep}$.
Let $\ka=\ba^+,$ 
the least cardinal bigger than $\ba$. 
Consider a countable elementary submodel $M$ of 
$\rL_\ka$ containing $\ba$, and let 
$H:M\onto \rL_\la$ be the Mostowski collapse. 
Let $\bbo=H(\ba)$; then $\bbo<\la$. 
As the construction of the sets $\rA_\al$ 
is obviously absolute for $\rL$, we have 
$\rA_{\bbo}\mo\text{\Sep}$ as well, and then 
$\rA_{\bbo}= \rA_{\bbo+1}$, as required.
\epf

The following theorem is essentially Lemma~2.2 
in \cite{marekG}.  

\bte
[$\ZF$]
\lam{b0pa} 
$\rA=\rA_{\bbo}$ satisfies\/ $\pad$ 
with the choice schema\/ $\ACw$.
\ete

To sketch a proof of this profound result, 
we need to 
have a look at the ramified analytic hierarchy   
from a somewhat different angle. 
This involves a \lap{shift} of 
G\"odel's hierarchy and ensuing classification of 
ordinals:
\bit
\item%
Let $\rM_\al=\rL_{\om+\al}$ for all $\al$. 
In particular, 
$\rM_0=\rL_\om=$ all hereditarily finite sets, 
but still, similarly to \eqref{LLdef}, 
$\rM_{\al+1}=\Def{\,\rM_\al}$, $\kaz\al$, and 
the union is taken at limit steps. 
(See \eg\ note 2 on p.\:499 in \cite{BPrah} or 
Section~5 in \cite{JSrah} where    
``$\rL_0=$ hereditarily finite sets'' 
is defined outright.) 
Needless to say that $\rM_\al=\rL_\al$ for all 
$\al\ge\om^2$.

\item%
An ordinal $\al$ is an \rit{index} if 
$(\rM_{\al+1}\bez\rM_\al)\cap\pws\om\ne\pu$,  
\cite{BPrah,LPrah,marekG}. 
\eit

\vyk{
\item[$-$]%
$\al$ \rit{starts a gap} \cite{marekG} if 
it is a gap and  
$\kaz\ga<\al,(\rL_{\al}\bez\rL_\ga)\cap\pws\om\ne\pu.$  
}

We'll refer to the following result, 
established in \cite{BPrah}, Theorems 1 and 9 
by a complex mixture 
of set theoretic (constructibility) 
and recursion theoretic methods. 
A set $E\sq\om\ti\om$ is \rit{a code\/} 
(or \rit{arithmetical copy\/}, as in 
\cite{BPrah,LPrah}) 
of $\rM_\al$ 
if it  
is isomorphic to ${\in}\res{\rM_\al}$ 
via a bijection of $\obl{E}$ onto $\rM_\al$. 

\bpro
\label{BP}
\ben
\renu
\itlb{BP2}%
If\/ $\al\le{\bbo+1}$ then\/ 
$\rA_\al=\rM_{\al}\cap\pws\om$. 

\itlb{BP1}%
If\/ $\ba$ is an index then there is a  
code of\/ 
$\rM_\ba$ in\/ $\rM_{\ba+1}$. 
\een
\epro 
\bpf[sketch]
\ref{BP1} 
Suppose that $\ba$ is limit. 
Argue as in Section~\ref{fura3} 
with $B=\om$ and 
$\rM_\ba=\rL_{\om+\ba}$ instead of $\rL$, so that 
$Y=\ens{k\in \om}{\vpi^{\rM_\ba}(k)}\nin \rM_\ba$.
In the notation of Section~\ref{fura3}, we still have 
$N=\rM_\la$ for a limit $\la$. 
Note that $\la<\ba$ is impossible since 
$Y\in \rM_{\la+1}\bez\rM_\ba$. 
And $\la>\ba$ is impossible as well since $N$ is 
the transitive collapse of $M\sq\rM_\ba$. 

Thus $\la=\ba$, and hence $\rM_\ba$ is 
$\in$-isomorphic to $M$. 

On the other hand, $M\in\rM_{\ba+1}$ as a definable 
subset of $\rM_\ba$. 
Moreover, the inductive construction of $M$ as the 
closude of $\om$ under a finite list of functions 
definable over $\rM_\ba$, can be represented as a 
construction of a relation $E\sq\om\ti\om$, still 
definable over $\rM_\ba$, and such that 
$\stk\om E$ is isomorphic to $\stk M{\in}$, hence 
to $\stk{\rM_\ba}{\in}$ by the above. 

In other words, $E\in\rM_{\ba+1}$ is a code of 
$\rM_\ba$, as required. 

If $\ba=\nu+k$, where $\nu$ is limit and $1\le k<\om$ 
then we have to go back to Section~\ref{fura3} and, 
using $\sg$, define a closed formula $\sg_k$ 
by induction on $k$, such 
that, for any transitive set $M$, 
$(\sg_k)^M$ holds iff $M=\rL_{\nu+k}$ for some 
limit ordinal $\nu$. 
Namely, put $\sg_0:=\sg$ as in Section~\ref{fura3}, 
then let $\sg_{k+1}$ say: 
``there is a transitive set $X$ with $(\sg_k)^X$ and 
(all sets) = $\Def X$''.

Then go through the arguments in 
the limit case, {\sl mutatis mutandis\/}. 

\ref{BP2} 
This claim goes by induction, using \ref{BP1} 
as the key argument. 
See \cite{BPrah} for details.
\vyk{%
The key step is as follows. 
Assume that $\rA_\al=\rM_{\al}\cap\pws\om$, 
and prove 
$\rA_{\al+1}=\rM_{\al+1}\cap\pws\om$. 
First of all, we immediately have 
$\rA_{\al+1}\sq\rM_{\al+1}\cap\pws\om$
by \eqref{LLdef} and \eqref{Adef}. 
It remains to show the inverse, \ie, if 
$x\sq\om$ belongs to $\rM_{\al+1}$ then 
$x\in \rA_{\al+1}$. 
By \ref{BP2}, 
consider a code $E\in\rM_{\ba+1}$, $E\sq\om^2$, 
of $\rM_\ba$. 
Then any $x\sq\om$ in $\rM_{\al+1}$ is arithmetical 
in $E$ (as it is definable over $\rM_\al$). 
}
\epF{Proposition}

\vyk{
The last equality is given in \cite{BPrah} 
as 
$\rA_\al=\rM_{\al}\cap\pws\om$, 
where sets $\rM_\al$ 
(modified constructible hierarchy)
are defined similarly to \eqref{LLdef} 
but starting  with 
}

\vyk{
\bpro
[\cite{BPrah}, Corollary 4]
\lam{3c4}
$\bbo$ is the least non-index.\qed
\epro 
}

\vyk{
\bpro
[\cite{marekG}, part of Theorem~2.7]
\lam{21t27}
If\/ $\al$ starts a gap then\/ 
$\rL_\al\mo\zfcm+\coun$.\qed
\epro 
}

\bpf
[Theorem~\ref{b0pa}]
The equality $\rA_{\bbo}= \rA_{\bbo+1}$ implies 
{\axf Comprehension} in $\rA_{\bbo}$. 
The proof of $\ACw$ 
takes more effort. 
We claim that:
\ben
\Renu
\itlb{172i}\msur%
{\em
$\bbo$ is not an index, whereas 
each\/ $\al<\bbo$ is an index}\/; 

\vyk{
\itlb{172ii}%
By \ref{172i} and Proposition~\ref{BP}\ref{BP1}, 
{\em
if\/ $\al<\bbo$ then 
there is a code\/ $x\in\rL_{\al+1}$ 
of\/ $\rL_{\al}$, 
and $\al$ is countable in\/ $\rL_{\al+1}$}.
}

\itlb{172iii}\msur%
{\em $\bbo$ is a limit ordinal\/} --- 
Lemma 2.5 in \cite{marekG}. 
\een

To prove \ref{172i}, 
note that, by the choice of $\bbo$ and 
Proposition~\ref{BP}\ref{BP2}, 
$\bbo$ is is not an index since  
$(\rM_{\bbo+1}\bez\rM_\bbo)\cap\pws\om=  
(\rA_{\bbo+1}\bez\rA_\bbo)\cap\pws\om=\pu$, 
whereas every $\al<\bbo$ is an index by 
similar reasons. 

To verify \ref{172iii}, suppose to the contrary 
that $\bbo=\al+1$. 
By \ref{172i} and 
Proposition~\ref{BP}\ref{BP1}, 
there is a code $x\sq\om$ of $\rM_\al$ 
in\/ $\rM_{\bbo}$, hence, in $\rA_{\bbo}$ 
by Proposition~\ref{BP}\ref{BP2}. 
In particular, $x$ codes all sets in 
$\rM_\al\cap\pws\om$. 
Therefore we can extract a part $y\sq\om$ 
of $x$, which codes all those sets so that 
\blin
\bul
{17*}
{%
\rM_\al\cap\pws\om=\ens{(y)_n}{n<\om}\ 
, 
}
\elin
\text{(see Section~\ref{prePA} on $(x)_n$)},
and in addition $y$ is arithmetical in $x$. 

Then $y\in\rA_\bbo$ by Lemma~\ref{bhp}. 
But each $z\in \rA_\bbo$ is 
arithmetical in $y$ by \eqref{17*}. 
This is a contradiction since 
$\rA_{\bbo}\mo\padm$ by Lemma~\ref{eb0}.

Now, coming to $\ACw$, we are going to prove that 
\blin
\bul
{17ac}
{\kaz n\,\sus x\,\Phi(n,x)\imp
\sus y\,\kaz n\,\Phi(n,(y)_n)}
\elin
holds in $\rA_\bbo$, where $\Phi$ is a 
$\pad$ formula possibly with parameters 
in $\rA_\bbo$. 

We make use of the \weo\ $\lcl\bbo$ of $\rM_{\bbo}$, 
definable over $\rM_{\bbo}$. 
Such an ordering 
exists by Lemma~\ref{13c3} since $\bbo$ 
is limit by \ref{172iii}. 
Assuming that the left-hand side of \eqref{17ac} 
holds in $\rA_\bbo$,  we let $x_n$ be 
the $\lcl\bbo$-least element 
$x\in \rA_\bbo=\rM_\bbo\cap\pws\om$ satisfying 
$\rA_\bbo\mo \Phi(n,x)$. 

The set $y=\ens{\asko nj}{j\in x_n}$ is then 
definable over $\rM_{\bbo}$, hence 
$y\in\Def{\,\rM_{\bbo}}=\rM_{\bbo+1}$. 
We conclude that  $y\in\rA_{\bbo+1}$ by 
Proposition~\ref{BP}\ref{BP2}. 
Finally $y\in\rA_{\bbo}$, because 
$\rA_{\bbo}=\rA_{\bbo+1}$ by the choice of $\bbo$. 
Thus $y$ witnesses the right-hand side of 
\eqref{17ac} since $(y)_n=x_n$ by construction. 
\epf

\vyk{

See 
\cite{%
JSrah,
welchRAH%
}
on some studies of the structure 
of ramified analytic sets.

}

The construction of the 
ramified analytical hierarchy 
is purely analytical and can be described by 
suitable 
$\lpad$ formulas. 
In principle, the proof of Theorem~\ref{b0pa} 
remains valid in $\zmc$ {\sl mutatis mutandis\/}. 
For instance, as $\omi$ may not exist 
in $\zmc$, the case $\bbo=\Ord$ has to be taken 
care of. 
Let
\blin
\bul
{anal1}
{\bbo=\left\{
\bay{lcl}
\text{the least }\ba
\text{ with }\rA_\ba=\rA_{\ba+1}&-&
\text{if such ordinals $\ba$ exist},\\[0.6ex]
\Ord, \text{ the class of all ordinals}&-&
\text{otherwise},
\eay
\right.}
\elin
so that $\rA=\bigcup_{\al\in\bbo}\rA_\al$ in 
both cases. 
It can be an interesting problem though to 
maintain the construction and the proof of 
Theorem~\ref{b0pa} entirely by analytical means 
on the base of $\padm$, thereby giving a pure 
analytical proof of the ensuing equiconsistency 
of $\padm$ and $\pad$.

\parf{Conclusions and problems}
\las{rp}

In this study, the methods of second-order arithmetic    
and set theory were employed to giving a full and 
self-contained proof of Theorem~\ref{1t} on 
the formal equiconsistency 
of such theories as second-order arithmetic $\padm$ 
and Zermelo--Fraenkel $\zfcm$ without the Power Set 
axiom (Theorem~\ref{1t}).  

\vyk{
This essentially strengthens and extends our earlier 
results in \cite{kl56}, 
by $\rV=\rL[a]$ and the minimality claim. 
In addition, we established (Theorem~\ref{1t}) 
the ensuing consistency result on the basis of second 
order Peano arithmetic $\pad$, instead of the much stronger 
theory $\ZFC$ typically assumed as a premise 
in independence results obtained by the forcing method. 
This is a new result and a valuable improvement upon 
much of known independence results in modern set theory.  

The technique developed in this paper may lead 
to further progress in studies of different aspects 
of the projective hierarchy.
We hope that this study will contribute 
to the following crucial problem by S.\,D.\,Friedman,    
see {\cite[P. 209]{sdffs} and  
\cite[P. 602]{sdf_ccf}}: 
\rit{find a model of\/ $\ZFC$, for a given\/ $\nn$, 
in which all\/ 
$\fs1\nn$ sets of reals are Lebesgue measurable 
and have the Baire and perfect set properties, 
and in the same time there exists a\/ $\fd1{\nn+1}$ 
\weo\ of the reals}. 
}

\vyk{
From our study, it is concluded that the 
technique of 
\rit{transitive models of bounded Separation in\/ $\zfcm$\/}, 
as in Section~\ref{2prel},   
will lead to similar consistency and independence 
results, related to second 
order Peano arithmetic $\pad$ and 
similar to our Theorem~\ref{1t}, 
on the basis of the consistency of $\pad$ itself.
}

The following problems arise from our study.

\begin{Problem}
\rm
\lam{pro1}
Regarding the axiom \trsu\  
(Transitive superset, Definition~\ref{zmad}), is 
it really independent of the rest of $\zmc$ axioms? 
\end{Problem}

\vyk{
\begin{Problem}
\rm
\lam{pro2}
The nature of the set or class $\rla$ of 
Section~\ref{fura} remains not fully clear. 
In $\ZF$, 
obviously $\xom\sq\omi$ is an ordinal, 
and accordingly 
$\rla$ is definitely a set satisfyng  
$\rla\sq\rL_{\omi}$. 
What about the possibility of the equality and 
of the strict inclusion $\sneq$ in this relation?
\end{Problem}

\begin{Problem}
\rm
\lam{pro3}
Consider, in $\zmc$, the structure  
$\stk\om{\pws\om\cap\rL}$ 
of all constructible subsets of $\om$. 
Does it satisfy $\padm$, or even stronger, $\pad$?
\end{Problem}
}

\begin{Problem}
\rm
\lam{pro1+}
Find a purely analytical proof of Theorem~\ref{b0pa} 
in $\padm$ that does not involve 
$\bI$ of Definition~\ref{bis}, 
or any similar derived set-theoretic structure, 
explicitly or implicitly.  
\end{Problem}

We expect that the methods and results of this 
paper can be used to strengthen and further 
develop Cohen's set-theoretic forcing method  
in its recent applications to theories 
$\zfcm$ and $\pad$ in \cite{kl75}. 
The technique of definable generic forcing notions 
has been recently applied for some definability 
problems in modern set theory, including 
the following applications$:$ 
\bit
\item[$-$]
a model of $\ZFC$ in \cite{kl36}, in which 
minimal collapse functions $\om\onto\omil$ 
first appear at a given projective level$;$  

\item[$-$]
a model of $\ZFC$ in \cite{kl49}, 
in which the Separation principle fails 
for a given projective class 
$\fs1\nn$, $\nn\ge3$$;$

\item[$-$]
a model of $\ZFC$ in \cite{kl29}, 
in which the full basis theorem 
holds in the absence of analytically definable 
\weo s of the reals$;$.

\item[$-$]
a model of $\ZFC$ in \cite{kl65}, 
in which the Separation principle 
holds for a given effective class 
$\is1\nn$, $\nn\ge3$. 
\eit
It is a common problem related to all these results to 
establish their $\pad$-consistency 
versions similar to Theorem~\ref{1t}.



\back
The authors are thankful to Ali Enayat, Gunter Fuchs,  
Victoria Gitman, and Kameryn Williams, 
for their enlightening comments that made it possible to 
accomplish this research. 
%
\vyk{
The authors are thankful to 
the anonymous referee for their 
comments and suggestions, which significantly 
contributed to improving the quality of the publication.
}
\eack



\renek{\refname}
{{\large\bf References}%
\addcontentsline{toc}{subsection}{References}}

{\small

\bibliographystyle{plain} 

\bibliography{80f.bib,80fkl.bib,80fmar.bib}

}

\small

\renek{\indexname}
{{\large\bf Index}%
\addcontentsline{toc}{subsection}{Index}}

\input{80find.sty}


\end{document}